\documentclass{amsart}
\usepackage[mathscr]{eucal}
\usepackage{amssymb, amsmath,array, amscd}
\usepackage[OT2,T1]{fontenc}

\newtheorem{thm}{Theorem}
\newtheorem{thma}{Theorem}[section]

\newtheorem{prop}{Proposition}

\newtheorem{lemma}{Lemma}[section]
\newtheorem{claim}{Claim}[section]

\newtheorem{rem}{Remark}[section]

\newtheorem{Cor}{Corollary}

\newtheorem{prob}{Problem}

\newtheorem{ex}{Example}

\def\wt{\widetilde}
\def\dd{\mathcal{D}}
\def\wh{\widehat}

\def\pr{{\,\prime}}
\def\Ga{\Gamma}

\def\p{\mathbb{P}}

\def\Te{\Theta}
\def\te{\theta}
\def\lim{\underset{z\to 0}{lim}\,}
\def\pa{\partial}

\def\sub{\subset}

\def\emp{\emptyset}
\def\ov{\overline}
\def\om{\omega}
\def\Om{\Omega}
\def\fa{\mathcal{F}}
\def\a{\alpha}
\def\be{\beta}

\def\C{\mathbb{C}}

\def\{{\lbrace}
\def\}{\rbrace}

\def\la{\lambda}

\def\L{\mathcal{L}}
\def\g{\gamma}

\def\N{\mathbb{N}}
\def\Si{\Sigma}
\def\te{\theta}

\def\G{\mathcal G}

\def\Z{\mathbb{Z}}
\def\O{\mathcal{O}}
\def\*{\star}

\def\Q{\mathbb{Q}}

\begin{document}

\title{Codimension two holomorphic foliations.}

\author{D. Cerveau \& A.  Lins Neto}

\begin{abstract}
This paper is devoted to the study of codimension two holomorphic
foliations and distributions. We prove the stability of complete intersection of codimension two distributions and foliations
in the local case. Converserly we show the existence of codimension two foliations
which are not contained in any codimension one foliation. We study problems related
to the singular locus and we classify homogeneous foliations of small
degree.
\end{abstract}



\maketitle

\tableofcontents

\section{Introduction}
There are many works devoted to the study of codimension one holomorphic foliations on complex manifolds. The local theory is well understood in small dimensions (2 and 3), with results concerning reduction of singularities ([Se], and [Cn-Ce]) and applications to unfolding theory, topological classification ([M-M-1] and [M-M-2]), spaces of moduli ([Ma]), existence and construction of invariant hypersurfaces ([C-M]), first integrals ([Ma-Mo] and [Ce-Ma]), among other topics.

In the global case, there is an intensive activity concerning the description of the "irreducible components" of the space of codimension one holomorphic foliations on a compact complex manifold ([Ce-LN], [CA-1] and [Ce-LN-CA-G]). One of the most popular challenges is to know if any codimension one foliation on $\p^n$, $n\ge3$, is either the meromorphic pull-back of a foliation on a complex surface, or has a "geometric" transverse structure ([Ce-LN-Lo-Pe-Tou] and [Ce-LN-2]).

In the present work, we focus our attention on singular foliations and distributions of codimension $q$, $q\ge2$, with special emphasis in the case $q=2$. Local and global results are obtained. For example, a way to construct a singular codimension two distribution is to intersect two singular codimension one distributions. In the local case we prove in theorem \ref{t:2} the "stability" of such construction, under natural assumptions. As a consequence, using Malgrange's singular Frobenius theorem, we show the persistence of first integrals (corollary \ref{c:1}).Conversely, we prove the existence of codimension two foliations which are not "contained" in any codimension one foliation; this fact is proved in the local contest and on rational manifolds (see proposition \ref{p:3}, corollary \ref{c:2} and remark \ref{r:25}).

Next the following problem is studied: is there a germ at $0\in\C^4$ of codimension two foliation with an isolated singularity at $0$?
Indeed, there are examples of holomorphic codimension two distributions on $\C^4$ with an isolated singularity at $0\in\C^4$. An example of this type was given in [K-N] in the contest of vector bundles on $\p^3$: it is defined by a homogeneous 2-form on $\C^4$, that is a 2-form with coefficients homogeneous of the same degree. In contrast with this example we prove in theorem \ref{t:5} that a codimension two foliation on $\C^4$, defined by a homogeneous 2-form with a singularity at $0\in\C^4$, has always a straight line in its singular set; in other words the singular set has dimension $\ge1$.

Finally, we describe with some details homogeneous foliation of small degree. That description is related to the classification of codimension one foliations of degree $\le2$ on $\p^n$, $n\ge3$.

\section{Definitions and some results}

\subsection{Local definitions}\label{se:11}
A holomorphic singular distribution of codimension $q$ (or dimension $n-q$) on an Stein open set $U\sub \C^n$, $0<q<n$, can be defined by a holomorphic $q$-form $\eta$ which is locally decomposable outside the singular set $Sing(\eta):=\{z\in U\,|\,\eta_z=0\}$, in the sense that any $z_o\in U\setminus Sing(\eta)$ has an open neighborhood $V\sub U$ such that
\begin{equation}\label{eq:1}
\eta|_V=\om_1\wedge...\wedge\om_q\,\,,
\end{equation}
where $\om_1,...,\om_q\in\Om^1(V)$. It follows that we can define in $U\setminus Sing(\eta)$ a holomorphic distribution $\dd_\eta$ of codimension $q$ by
\[
\dd_\eta(p)=\{v\in T_pU\,|\,i_v\eta(p)=0\}\,\,.
\]
If $p\in U\setminus Sing(\eta)$ and $\om_1,...,\om_q$ are as in (\ref{eq:1}) then
\[
\dd_\eta(p)=\bigcap_{1\le j\le q}Ker(\om_j(p))\,\,.
\]
A q-form $\eta$ satisfying (\ref{eq:1}) is said to be locally decomposable.

A q-form $\eta$ which satisfies (\ref{eq:1}) is integrable if it satifies Frobenius' integrability condition:
\begin{equation}\label{eq:2}
d\om_j\wedge\eta=0\,\,,\,\,\forall\,\,j=1,...,q\,,\,\text{on the open set}\,\,U\,\,.
\end{equation}

If $\eta$ satisfies (\ref{eq:1}) and (\ref{eq:2}) then the distribution $\dd_\eta$ is integrable and so $\eta$ defines a holomorphic codimension $q$ foliation on $U\setminus Sing(\eta)$. This foliation will be denoted by $\fa_\eta$.

The leaves of $\fa$ are the immersed codimension $q$ submanifolds $\L\sub U\setminus Sing(\eta)$ for which the tangent space to $\L$ at $m\in U$ is $T_m\L:=\dd_\eta(m)$.

When $Sing(\eta)=\emp$ then the integrability condition (\ref{eq:2}) is equivalent to the existence of a holomorphic 1-form $\te$ such that
\begin{equation}\label{eq:3}
d\eta=\eta\wedge\te\,\,.
\end{equation}
When $Sing(\eta)\ne\emp$ then (\ref{eq:3}) is only true locally in $U\setminus Sing(\eta)$, unless we allow $\te$ to be meromorphic.

\begin{ex}
Complete intersection. {\rm Let $\fa_1,...,\fa_q$ be $q$ foliations of codimension one defined by integrable 1-forms $\om_1,...,\om_q\in\Om^1(U)$, $1\le q<n$ ($\om_k\wedge d\om_k=0$, $\forall k$) such that $\eta:=\om_1\wedge...\wedge\om_q(m)\not\equiv0$. The foliation complete intersection $\fa=\fa_1\cap...\cap\fa_q$ is associated to the $q$ 1-forms $\om_i$, $1\le i\le q$. The leaves of $\fa$ are the connected components of the intersection of the leaves $\L_k$ of $\fa_k$, $1\le k\le q$.}
\end{ex}

\begin{ex}
Foliations associated to a Lie algebra of vector fields or to an action of a Lie group. {\rm Let $\L$ be a Lie algebra of vector fields defined on an open Stein subset $U\sub\C^n$. Given $m\in U$ set
\[
d(m)=dim_\C\left<X(m)\,|\,X\in\L\right>\,\,.
\]
Let $d=max\{d(m)\,|\,m\in U\}$ and assume that $1\le d<n$. 
Since $U$ is connected, then the set $Z=\{m\in U\,|\,d(m)<d\}$ is a proper analytic subset of $U$, so that $V=U\setminus Z$ is open dense in $U$ and connected. In particular, $\L$ defines a dimension $d$ distribution on $V$
\[
\L(m)=\left<X(m)\,|\,X\in \L\right>\,\,,\,\,m\in V\,\,.
\]
The Lie algebra $\L$ defines a codimension $q=n-d$ foliation $\fa_\L$ on $V$.
The foliation $\fa_\L$ can be extended to $U$ as a singular foliation with singular set $Z$.

When $\L$ is associated to a group action $G\times U\to U$ we will say that $\fa_\L$ is associated to the action of $G$.}
\end{ex}

\begin{rem}
{\rm If $\eta\in\Om^k(U)$ is integrable and $cod(Sing(\eta))=1$ then we can write $\eta=h.\,\eta^\pr$, where $\eta^\pr\in \Om^k(U)$ and $cod(Sing(\Om^\pr))\ge2$. We would like to observe that $\eta^\pr$ is also integrable. The foliation $\fa_{\eta^\pr}$ can be considered as an 
"extension" of $\fa_\eta$.}
\end{rem}

\begin{rem}\label{r:dom}
{\rm If $\eta$ is an integrable q-form and $d\eta\not\equiv0$ then relation (\ref{eq:3}) implies that $d\eta$ is locally decomposable outside $Sing(\eta)$. Since $d\eta$ is closed it is integrable and defines a singular foliation of codimension $q+1$. Relation (\ref{eq:3}) implies also that any leaf of $d\eta$ is $\eta$-invariant, in the sense that, either it is contained in a leaf of $\eta$, or in $Sing(\eta)$.}
\end{rem}

\begin{ex}\label{ex:3}
{\rm Let $\fa_1$ and $\fa_2$ be the codimension one foliations of $\C^n$ defined by the 1-forms $\eta_j=x_1...x_n.\,\om_j$, $j=1,2$, where $\om_1$ and $\om_2$ are the logarithmic closed forms
\[
\om_1=\sum_{j=1}^n\la_j\frac{dx_j}{x_j}\,\,,\,\,\om_2=\sum_{j=1}^n\mu_j\frac{dx_j}{x_j}\,\,,\,\,\la_j\,,\,\mu_j\in\C^*\,,\,1\le j\le n\,\,.
\]
We assume that $\om_1$ and $\om_2$ are not colinear, which is equivalent to the non-colinearity of the vectors $\la=(\la_1,...,\la_n)$ and $\mu=(\mu_1,...,\mu_n)$. The intersection $\fa_1\cap\fa_2$ is defined outside $\Si:=\bigcup_j\{x_j=0\}$ by 
\[
\om_1\wedge\om_2=\sum_{i<j}(\la_i\mu_j-\la_j\mu_i)\,\frac{dx_i}{x_i}\wedge\frac{dx_j}{x_j}\,\,.
\]
Note that $\om_1\wedge\om_2(m)\ne0$, for all $m\in \C^n\setminus\Si$.
Moreover, the divisor of poles of $\om_1\wedge\om_2$ is $x_1...x_n$ and the form
\[
\eta:=x_1...x_n.\,\om_1\wedge\om_2=\sum_{i<j}(\la_i\,\mu_j-\la_j\,\mu_i)\,x_1...\wh{x_i}...\wh{x_j}...x_n\,dx_i\wedge dx_j
\]
is holomorphic on $\C^n$, so that $\eta$ defines the codimension two foliation $\fa_1\cap\fa_2:=\fa$ on $\C^n$.
By convention, $\wh{x_i}$ means the omission of the factor $x_i$ in the product.

Observe that the hyperplanes $(x_j=0)$, $1\le j\le n$, are $\fa$-invariant. For instance, when $j=1$ we have $\eta|_{(x_1=0)}=dx_1\wedge\eta_1$, where
\[
\eta_1=\sum_{j>1}(\la_1\,\mu_j-\la_j\,\mu_1)\,x_2...\wh{x_j}...x_n\,dx_j\,\,,
\]
and $\eta_1\not\equiv0$, because otherwise $\la$ and $\mu$ would be colinear. In particular, $Sing(\eta)$ has codimension $\ge2$ and is contained in $\Si$.
In fact, if $\la_i\,\mu_j-\la_j\,\mu_i\ne0$ for all $i<j$ then
\[
Sing(\eta)=\bigcup_{i<j<k}(x_i=x_j=x_k=0)
\]
has codimension three.

Another fact, is that $\eta$ is not decomposable, whereas $\frac{1}{x_1...x_n}\,\eta=\om_1\wedge\om_2$.}
\end{ex}

\subsection{Globalization}\label{se:12}
Let $M$ be a complex manifold of dimension $n\ge2$. A holomorphic (singular) foliation $\fa$ of codimension $q$,$1\le q<n$, on $M$ is defined by a covering $(U_j)_{j\in J}$ of $M$ by Stein open sets, a collection of integrable $q$-forms $(\eta_j)_{j\in J}$, with $\eta_j\in \Om^q(U_j)$ and $cod(Sing(\eta_j))\ge2$, and a collection $(g_{ij})_{U_i\cap U_j\ne\emp}$ with $g_{ij}\in\O^*(U_i\cap U_j)$, satisfying the glueing condition $\eta_i=g_{ij}\,\eta_j$ on $U_i\cap U_j\ne\emp$. If $U_i\cap U_j\ne\emp$ then the glueing condition implies that the leaves of $\fa_{\eta_i}|_{U_i\cap U_j}$ coincides with the leaves of $\fa_{\eta_j}|_{U_i\cap U_j}$ and that  
\begin{equation}\label{eq:4}
Sing(\eta_i)|_{U_i\cap U_j}=Sing(\eta_j)|_{U_i\cap U_j}\,\,.
\end{equation}
Relation (\ref{eq:4}) implies that
\[
Sing(\fa):=\bigcup_{j\in J}Sing(\eta_j)
\]
is an analytic subset of $M$ of codimension $\ge2$.
Observe that this defines a non-singular foliation of codimension $q$ on the set $M\setminus Sing(\fa)$.

\subsection{Complete intersection of foliations}\label{se:13}
Let $\fa_1,...,\fa_q$ be codimension one foliations on some polydisc $U\sub\C^n$, defined by integrable 1-forms $\om_1,...,\om_q$ such that $\eta^\pr:=\om_1\wedge...\wedge\om_q\not\equiv0$ (generically independent). Note that there exists $f\in\O(U)$ and $\eta\in\Om^q(U)$ such that $\eta^\pr=f.\,\eta$ and $cod(Sing(\eta))\ge2$. Clearly, $\eta$ is integrable.

By definition, the foliation $\fa_\eta$, defined by $\eta$, is the {\it topological intersection} of the foliations $\fa_1,...,\fa_q$. When $cod(Sing(\eta^\pr))\ge2$ (or equivalently $f\in\O^*(U)$) we will say that $\fa$ is a {\it complete intersection}. 
Clearly, these two definitions can be germified or globalized (as in \S \ref{se:12}).

\begin{rem}
{\rm We will see in \S $\,$\ref{se:17} that there are germs of codimension two foliations $\fa$ that are not topological intersections of two codimension one foliations. In this case, $\fa$ is defined by a germ of integrable 2-form, say $\eta$, which is meromorphically decomposable, $\eta=\a\wedge\be$, but for any such decomposition, neither $\a$ nor $\be$ is integrable.}
\end{rem}

A direct consequence of a theorem due to Malgrange (cf. [M-1] and [M-2]) is the following:

\begin{thma}\label{t:1}
Let $\eta$ be a germ at $0\in\C^n$ of integrable q-form holomorphically decomposable, $\eta=\a_1\wedge...\wedge\a_q$. If $cod(Sing(\eta))\ge3$ then $\fa_\eta$, the foliation defined by $\eta$, is a complete intersection. More precisely, there are $f_1,...,f_q\in\O_n$ and $h\in\O^*_n$ such that
\[
\eta=h.\,df_1\wedge...\wedge df_q\,\,.
\]
\end{thma}

On the other hand, as we have seen in example \ref{ex:3} of \S$\,$ \ref{se:11}, for generic values of $\la_i$ and $\mu_i$, the form $\eta=\sum_{i<j}(\la_i\mu_j-\la_j\mu_i)\,x_1...\wh{x_i}...\wh{x_j}...x_n\,dx_i\wedge dx_j$ is integrable and satisfies $cod(Sing(\eta))\ge3$, but $\fa_\eta$ is not a complete intersection. Therefore, the hypothesis of holomorphic decomposability in Malgrange's theorem is crucial.

\subsection{Special case: $q=2$ and $n=4$}\label{se:14}
Let $\eta$ be a 2-form on a polydisc $U\sub \C^4$. We fix coordinates $x_1,x_2,x_3,x_4$ and a non-vanishing 4-form, for instance $\nu=dx_1\wedge...\wedge dx_4$. The 3-form $d\eta$ can be written as 
\[
d\eta=i_X\,\nu\,\,,
\]
where $X$ is a holomorphic vector field on $U$, called the {\it rotational} of $\eta$: $X=rot(\eta)$. The foliation associated to the rotational is independent of the choice of $\nu$. On the other hand, if we change $\eta$ by $h\,\eta$, where $h\in\O^*(U)$, then $rot(h\,\eta)$ and $rot(\eta)$ are not in general colinear. Although this notion is not intrinsic, it is convenient (see [LN]). For instance, in the two following propositions, the first valid in any dimension. 

\begin{prop}\label{p:1}
Let $\eta\in\Om^2(U)$, a 2-form on an open set $U$ of $\C^n$. Then $\eta$ is meromorphically decomposable if, and only if, $\eta^2=\eta\wedge\eta=0$.
\end{prop}

{\it Proof.}
If $\eta=\om_1\wedge\om_2$, $\om_1$ and $\om_2$ meromorphic, then clearly $\eta^2=0$. Conversely, if $\eta^2=0$ (and $\eta\not\equiv0$) then we choose two vector fields $Z_1$ and $Z_2$ such that $i_{Z_1}i_{Z_2}\,\eta\not\equiv0$. From $\eta^2=0$ we get
\[
i_{Z_1}\eta\wedge\eta=0\,\,\implies\,\,(i_{Z_1}i_{Z_2}\eta)\,\eta-i_{Z_1}\eta\wedge i_{Z_2}\eta=0\,\,\implies\,\,\eta=\om_1\wedge\om_2\,\,,
\]
where $\om_1=\frac{1}{i_{Z_1}i_{Z_2}\eta}\,i_{Z_1}\eta$ and $\om_2=i_{Z_2}\eta$.
\qed

\begin{prop}\label{p:2}
Let $U$ be a domain of $\C^4$ and $\eta\in\Om^2(U)$, $\eta\not\equiv0$. If $rot(\eta)\not\equiv0$ then $\eta$ is integrable if, and only if, $i_{rot(\eta)}\eta=0$.
\end{prop}

{\it Proof.}
Let us denote $Y:=rot(\eta)$; $d\eta=i_Y\nu$, $\nu=dz_1\wedge...\wedge dz_4$.
Since $\eta\in\Om^2(U)$, we can write $\eta^2=F.\,\nu$, where $F\in\O(U)$.
If $i_Y\eta=0$ then
\[
i_Y\eta^2=0\,\,\implies\,\,F.\,i_Y\nu=0\,\,\overset{Y\ne0}\implies\,\,F=0\,\,\implies\,\,\eta^2=0\,\,.
\] 
Therefore, by proposition \ref{p:1} we can write $\eta=\om_1\wedge\om_2$, where $\om_1$ and $\om_2$ are meromorphic.
We are going to prove that $d\om_1\wedge\eta=d\om_2\wedge\eta=0$ (which implies the integrability). 
From $i_Y\eta=0$ we have
\[
0=i_Y(\om_1\wedge\om_2)=i_Y(\om_1).\,\om_2-i_Y(\om_2).\,\om_1\,\,\implies\,\,i_Y\,\om_1=i_Y\,\om_2=0\,\,,
\]
because $\om_1$ and $\om_2$ are linearly independent in an open and dense set. On the other hand $\om_1\wedge\nu=0$ and so
\[
0=i_Y(\om_1\wedge\nu)=-\om_1\wedge i_Y\nu=-\om_1\wedge d\eta=-\om_1\wedge(d\om_1\wedge\om_2-\om_1\wedge d\om_2)\,\implies
\]
\[
d\om_1\wedge\eta=d\om_1\wedge\om_1\wedge\om_2=\om_1\wedge d\om_1\wedge\om_2=0\,\,.
\]
Similarly $d\om_2\wedge\eta=0$.
The converse is left to the reader.
\qed

\begin{rem}\label{r:13}
{\rm If $rot(\eta)=0$, i.e. if $\eta$ is closed, then $\eta$ is integrable if, and only if, $\eta^2=0$. The proof can be found in [LN].}
\end{rem}

\begin{ex}\label{e:4}
Codimension two distibutions and the generic quadric of $\,\p^5$. {\rm A 2-form $\eta$ on an open set $U\sub\C^4$ can be written as
\[
\eta=A\,dx_2\wedge dx_3+B\,dx_3\wedge dx_1+C\,dx_1\wedge dx_2+(E\,dx_1+F\,dx_2+G\,dx_3)\wedge dx_4\,\,.
\]
As the reder can check directly, the condition $\eta^2=0$ is equivalent to:
\[
A\,E+B\,F+C\,G=0
\]
Therefore, $\eta$ defines a singular codimension two distribution on $U$ if, and only if, the image of the map $\Phi=(A,B,C,D,E,F)\colon U\to\C^6$ is contained in the quadric $Q=(z_0\,z_3+z_1\,z_4+z_2\,z_5=0)$. When the components of $\Phi$ are homogeneous polynomials of the same degree then $\Phi$ defines a rational map $\phi\colon \p^3-\to \wh{Q}$, where $\wh{Q}\sub\p^5$ is the projection on $\p^5$ of the quadric $Q$. The indetermination set of $\phi$ is of course the projection on $\p^3$ of the set $\Phi^{-1}(0)$. In \S \ref{se:22} we will study the homogeneous case with $\Phi^{-1}(0)=\{0\}$. This corresponds to the case in which the form $\eta$ has an isolated singularity at $0\in\C^4$.}
\end{ex}

\subsection{Singularities and the rotational}\label{se:15}
Let $\eta$ be a germ at $0\in\C^4$ of integrable 2-form with a singularity at $0$. We will examine two cases:

\subsubsection{The Kupka-Reeb phenomenon}\label{s:151}
When $rot(\eta)(0)\ne0$ we can find local coordinates $x=(x_1,x_2,x_3,x_4)$ around $0\in\C^4$  such that $rot(\eta)=\frac{\pa}{\pa x_4}$. In these coordinates the form $\eta$ does not depends on the variable $x_4$ and of $dx_4$. More precisely, we can write $\eta=i_Z\,dx_1\wedge dx_2\wedge dx_3$, where $Z$ is a germ of vector field as below
\begin{equation}
Z=\sum_{j=1}^3A_j(x_1,x_2,x_3)\,\frac{\pa}{\pa x_j}\,\,.
\end{equation}
In particular, the foliation $\fa_\eta$ can be interpreted as the pull-back by the projection $x\mapsto (x_1,x_2,x_3)$ of the germ of foliation on $(\C^3,0)$ defined by the vector field $Z$ (cf. [Me]).

\subsubsection{The case in which $rot(\eta)$ has an isolated singularity}
When $rot(\eta)(0)=0$ and $0$ is an isolated singularity of $rot(\eta)$ we can apply the division theorem of De Rham-Saito [DR] as follows:
the integrability condition $i_{rot(\eta)}\eta=0$ implies that there exists a germ of holomorphic vector field $S$ such that
\begin{equation}\label{eq:6}
\eta=i_S\,i_{rot(\eta)}\,dx_1\wedge dx_2\wedge dx_3\wedge dx_4\,\,.
\end{equation}
In other words, the tangent bundle of $\fa$ decomposes globally outside $Sing(\fa_\eta)$: $T\fa=\C.\,S\oplus \C.\,rot(\eta)$.
Note that (\ref{eq:6}) implies that $L_S\eta=\eta$, where $L$ denotes the Lie derivative. 

Two sub-cases were studied in [LN]:

\begin{itemize}
\item[1$^{st}$.] When the linear part of $rot(\eta)$ at $0$ is non-degenerate. In this case, under generic conditions, it is possible to find germs of vector fields $X,Y$ generating the $\O_4$-module $T\fa$ and such that $[X,Y]=0$: the foliation is generated by a local action of $\C^2$.
\item[2$^{nd}$.] When the linear part of $rot(\eta)$ at $0$ is nilpotent. In this case, it is proved in [LN] that the eigenvalues of $DS(0)$ are rational positive. In particular, $S$ is holomorphically normalizable; $S=S_1+N$, where $S_1=DS(0)$ and $N$ is nilpotent and $[S_1,N]=0$.
It is proved also that there exist local coordinates $x=(x_1,x_2,x_3,x_4)$ in which the rotational satisfies $[S_1,rot(\eta)]=\la\,rot(\eta)$, where $\la\in\Q_+$. In other words, $S$ and $rot(\eta)$ generate a local action of the affine group of $\C$.
\end{itemize}

\subsection{Stability of complete intersections. Case of 2-forms}\label{se:16}
This section is devoted to the following result:
\begin{thm}\label{t:2}
Let $\eta_0$ be a germ at $0\in\C^n$ of holomorphic decomposable 2-form, $\eta_0=\a_0\wedge\be_0$, where $\a_0,\be_0\in \Om^1(\C^n,0)$.
Assume that $cod(Sing(\eta_0))\ge3$. Let
\[
\eta_s=\eta_0+s\,\eta_1+...=\sum_{j\ge0}s^j\,\eta_j
\]
be a holomorphic family of germs of 2-forms such that $\eta_s^2=0$. Then $\eta_s$ is decomposable, i.e. there exist holomorphic families of 1-forms $\a_s=\a_0+s\,\a_1+...$ and $\be_s=\be_0+s\,\be_1+...$ such that $\eta_s=\a_s\wedge\be_s$.
\end{thm}

{\it Proof.}
It is sufficient to prove that we can write $\eta_s=\a_s\wedge\be_s$, where $\a_s=\sum_{j\ge0}s^j\,\a_j$ and $\be_s=\sum_{j\ge0}s^j\,\be_j$ are formal power series.
In fact, the existence of a formal solution of the equation $\,\eta_s-\a_s\wedge\be_s=0\,$ implies the existence of a convergent solution by Artin's approximation theorem (cf. [Ar]).

Given a power series $\te_s=\sum_{s\ge0}s^j\,\te_j$ we set $j^k(\te_s)=\sum_{j=0}^ks^j\,\te_j$. 
We will prove by induction on $k\ge0$ that there exist germs $\a_s^k:=\sum_{j=0}^ks^j\,\a_j$ and $\be_s^k:=\sum_{j=0}^ks^j\,\be_j$ such that 
\begin{equation}\label{eq:7}
j^k\left(\eta_s-\a_s^k\wedge\be_s^k\right)=0\,\,.
\end{equation}
The first step of the induction is the hypothesis $j^{\,0}(\eta_s-\a_0\wedge\be_0)=0$.
Suppose by induction that there exist $\a_s^{\ell-1}$ and $\be_s^{\ell-1}$ satisfying (\ref{eq:7}) for $k=\ell-1\ge0$, and let us prove that there exist $\a_s^\ell$ and $\be_s^\ell$ satisfying (\ref{eq:7}), $j^{\ell-1}(\a_s^\ell)=\a_s^{\ell-1}$ and $j^{\ell-1}(\be_s^\ell)=\be_s^{\ell-1}$.
Observe first that
\[
j^\ell(\eta_s-\a_s^{\ell-1}\wedge\be_s^{\ell-1})=s^\ell\,\mu\,\,,
\]
where
\[
\mu=\eta_\ell-\sum_{i=1}^{\ell-1}\a_i\wedge\be_{\ell-i}\,\,.
\]
By Saito's division theorem the induction step can be reduced to the following (cf. [Sa]):
\begin{claim}\label{cl:11}
If $\mu$ is as above then
$\a_0\wedge\be_0\wedge\mu=0$.
\end{claim}

In fact, if claim \ref{cl:11} is true then, since $cod(Sing(\a_0\wedge\be_0))\ge3$, by Saito's theorem there are germs $\a_\ell$ and $\be_\ell$ such that
\[
\mu=\a_0\wedge\be_\ell+\a_\ell\wedge\be_0\,\,.
\]
Therefore, if we set $\a_s^\ell=\a_s^{\ell-1}+s^\ell\,\a_\ell$ and $\be_s^\ell=\be_s^{\ell-1}+s^\ell\,\be_\ell$ then
\[
j^\ell(\eta_s-\a_s^\ell\wedge\be_s^\ell)=0\,\,,
\]
as the reader can verify directly.
\vskip.1in

{\it Proof of claim \ref{cl:11}.}
Here we use the relation $\eta_s^2=0$, which implies
\[
0=\eta_s^2=\sum_rs^r\,\sum_{m+n=r}\eta_m\wedge\eta_n\,\,\implies\,\,\sum_{m+n=r}\eta_m\wedge\eta_n=0\,\,,\,\,\forall r\ge0\,\,.
\]

The induction hypothesis implies
\[
\sum_{j=0}^{\ell-1}s^j\,\eta_j=j^{\ell-1}(\eta_s)=j^{\ell-1}(\a_s^{\ell-1}\wedge\be_s^{\ell-1})\,\,\implies\,\,
\eta_r=\sum_{i=0}^r\a_i\wedge\be_{r-i}\,\,,\,\,0\le r\le \ell-1\,\,.
\]
Therefore,
\begin{equation}\label{eq:8}
0=\sum_{m+n=\ell}\eta_m\wedge\eta_n=2\,\eta_0\wedge\eta_\ell+\sum_{\underset{r+s\ge1}{\underset{i+j\ge1}{i+j+r+s=\ell}}}\a_i\wedge\be_j\wedge\a_r\wedge\be_s
\end{equation}

Since $\a_i\wedge\a_i=0$ and $\be_i\wedge\be_i=0$, $\forall i$, we can assume that the summation set in the sum of the right hand side of (\ref{eq:8}) is
\[
S=\{(i,j,r,s)\,|\,i+j+r+s=\ell\,,\,i+j\ge1\,,\,r+s\ge1\,\,,\,\,i\ne r\,\,,\,\,j\ne s\}\,\,.
\]
For simplicity of notation, given a subset $A\sub S$ we will set
\[
\sum_{(i,j,r,s)\in A}\a_i\wedge\be_j\wedge\a_r\wedge\be_s:=\sum(A)\,\,.
\]

The set $S$ can be decomposed into two disjoint subsets $S=S_1\cup S_2$, where $S_1=\{(i,j,r,s)\in S\,|\,i=s=0\,\,\text{or}\,\,j=r=0\}$ and $S_2=S\setminus S_1$.
Note that $S_2=\{(i,j,r,s)\in S|$ at most one of the indexes $i$, $j$, $r$ or $s$ is $=0\}$.

We assert that $\sum(S_2)=0$.
This assertion follows from the fact that there exists a permutation $\phi\colon S_2\to S_2$ with the following properties:
\begin{itemize} 
\item[(i).] $\phi(x)\ne x$ and $\phi(\phi(x))=x$, $\forall x\in S_2$.
\item[(ii).] If $x=(i,j,r,s)\in S_2$ and $\phi(x)=(i^\pr,j^\pr,r^\pr,s^\pr)$ then
\[
\a_i\wedge\be_j\wedge\a_r\wedge\be_s=-\a_{i^\pr}\wedge\be_{j^\pr}\wedge\a_{r^\pr}\wedge\be_{s^\pr}
\]
\end{itemize}

In fact, (i) and (ii) imply that $\phi$ induces a partition $S_2=S_3\cup S_4$, $\phi(S_3)=S_4$, $\phi(S_4)=S_3$, such that $\sum(S_3)=-\sum(S_4)$. Therefore, $\sum(S_2)=\sum(S_3)+\sum(S_4)=0$.

The permutation $\phi$ is constructed as follows: fix $(i,j,r,s)\in S_2$. We have three possibilities:
\begin{itemize}
\item[1.] $i,j,r,s\ge1$. In this case we set $\phi(i,j,r,s)=(r,j,i,s)$.
\item[2.] $i=0$ or $r=0$. In this case we set $\phi(i,j,r,s)=(i,s,r,j)$.
\item[3.] $j=0$ or $s=0$. In this case we set $\phi(i,j,r,s)=(r,j,i,s)$. 
\end{itemize}
We leave to the reader the verification that $\phi$ is well defined and satisfies properties (i) and (ii).

Finally we can write $S_1=S_5\cup S_6$ where $S_5=\{(i,0,0,\ell-i)\,|\,1\le i\le \ell-1\}$ and $S_6=\{(0,\ell-i,i,0)\,|\,1\le i\le \ell-1\}$.
Since
\[
\a_i\wedge\be_0\wedge\a_0\wedge\be_j=\a_0\wedge\be_j\wedge\a_i\wedge\be_0=-\a_0\wedge\be_0\wedge\a_i\wedge\be_j
\]
we get
\[
\sum(S_1)=-2\,\a_0\wedge\be_0\wedge\sum_{i=1}^{\ell-1}\a_i\wedge\be_{\ell-i}\,\,\overset{(\ref{eq:8})}\implies\,\,\a_0\wedge\be_0\wedge\mu=0\,\,\qed
\]
\vskip.1in

Theorem \ref{t:2} has the following consequence:
\begin{Cor}\label{c:1}
Let $\eta_0$ be a germ at $0\in\C^n$ of decomposable and integrable 2-form: $\eta_0=\a_0\wedge\be_0$, where $cod(Sing(\eta_0))\ge3$. Let $\eta_s$, $s\in(\C,0)$, be a holomorphic family of germs of integrable 2-forms such that $\eta_0=\eta_s|_{s=0}$. Then there exist holomorphic families of germs of functions
$f_s$, $g_s$ and $u_s$ of germs of holomorphic functions, $u_s\in\O^*_n$, such that
$\eta_s=u_s\,df_s\wedge dg_s$. In particular, the foliation associated to $\eta_s$ has two independent first integrals.
\end{Cor}

{\it Proof.} Theorem \ref{t:2} implies the existence of two holomorphic families of 1-forms $\a_s$ and $\be_s$
such that $\eta_s=\a_s\wedge\be_s$, $\forall s\in(\C,0)$. Consider on $(\C^n\times\C,(0,0))$ the pffafian system generated by the forms $\a_s$, $\be_s$ and $ds$. This system is integrable. Since $cod_{(\C^n,0)}(Sing(\a_0\wedge\be_0))\ge3$, by semi-continuity we have $cod_{(\C^n,0)}(Sing(\a_s\wedge\be_s))\ge3$. Since $\a_s\wedge\be_s$ does not contain terms with $ds$ we can conclude that
\[
cod_{(\C^n\times\C,0}(Sing(\a_s\wedge\be_s\wedge ds))\ge3\,\,.
\] 
Therefore, by Malgrange's theorem (see theorem \ref{t:1}) there exist $F,G\in\O_{n+1}$ and $U\in\O^*_{n+1}$ such that
\[
\a_s\wedge\be_s\wedge ds=U\,dF\wedge dG\wedge ds\,\,.
\]
Hence, we can take the families $f_s$, $g_s$ and $u_s$ as
\[
f_\tau:=F|_{s=\tau}\,,\,g_\tau:=G|_{s=\tau}\,\text{and}\,u_\tau:=U|_{s=\tau}\,\,\qed
\]

\subsection{Codimension two foliations not contained in a codimension one foliations}\label{se:17}
In the construction of the examples we will use a result due to X. Gomez-Mont and I. Luengo:
\begin{thma}\label{t:11}
There exists a polynomial vector field on $\C^3$ with an isolated singularity at $0\in\C^3$ and without germ of analytic invariant curve through $0$.
\end{thma}

A consequence of theorem \ref{t:11} is the following:

\begin{prop}\label{p:3}
Let $Z$ be a germ at $0\in\C^3$ of vector field with an isolated singularity at $0$ and without germ of invariant curve through $0$.
Then $Z$ cannot be tangent to a germ at $0\in\C^3$ of holomorphic codimension one foliation.
\end{prop}

{\it Proof.}
Suppose by contradiction that $Z$ is tangent to some a germ at $0\in\C^3$ of codimension one foliation $\fa$.  
Let $\om$ be a germ of integrable 1-form defining $\fa$ and with $cod(Sing(\om))\ge2$. 
We assert that $Sing(\om)=\{0\}$.

Note first that the tangency condition is equivalent to $i_Z\om=0$, which implies that $Sing(\om)$ is $Z$-invariant.
Therefore, $Sing(\om)$ cannot contain a germ of curve, for otherwise this curve would be $Z$-invariant. Hence, $Sing(\om)\sub\{0\}$ and we have two possibilities, either $Sing(\om)=\emp$, or $Sing(\om)=\{0\}$. On the other hand, if $Sing(\om)$ was empty then by integrability there exists a local chart $x=(x_1,...,x_n)$ such that $\om=u\,dx_1$, where $u(0)\ne0$. This implies that $i_Zdx_1=Z(x_1)=0$ and this contradicts the fact that $0$ is an isolated singularity of $Z$. Therefore, $Sing(\om)=\{0\}$ as asserted.

Let $\eta=i_Z\,dx_1\wedge dx_2\wedge dx_3$ and observe that the relation $i_Z\om=0$ is equivalent to $\om\wedge\eta=0$. Since $cod(Sing(\om))=3$ by De Rham's theorem [DR] there exists a germ of 1-form $\te$ such that $\eta=\om\wedge\te$.
However, a decomposable 2-form $\eta=\om\wedge\te$ with a singularity at $0\in\C^3$ vanishes necessarily on a curve through $0$. This contradicts the fact that $Sing(\eta)=Sing(Z)=\{0\}$.
\qed

\begin{Cor}\label{c:2}
For all $n\ge3$ there are germs at $0\in\C^n$ of holomorphic codimension two foliations which are "not contained" in any holomorphic foliation of codimension one: if a germ like this is defined by an integrable 2-form $\eta$ then there is no integrable 1-form $\om$ such that $\om\wedge\eta=0$.
\end{Cor}

{\it Proof.} If $n=3$ the corollary is a direct consequence of proposition 3: take $\eta=i_Z\,dx_1\wedge dx_2\wedge dx_3$, where $Z$ is like in theorem \ref{t:11}. If $n>3$ then let $\Pi\colon\C^n=\C^3\times\C^{n-3}\to\C^3$ be the projection $\Pi(x,y)=x$ and
\[
\eta=\Pi^*(i_Z\,dx_1\wedge dx_2\wedge dx_3)\,\,.
\]

Suppose by contradiction that there exists a germ of integrable 1-form $\om$ such that $\om\wedge\eta=0$. Note that in the coordinates $(x,y)\in\C^3\times\C^{n-3}$, $y=(y_1,...,y_{n-3})$ the definition of $\eta$ implies that $i_{\frac{\pa}{\pa y_j}}\eta=0$,
$\forall\, 1\le j\le n-3$.
Therefore,
\[
0=i_{\frac{\pa}{\pa y_j}}(\om\wedge\eta)=i_{\frac{\pa}{\pa y_j}}(\om).\,\eta\,\,\implies\,\,i_{\frac{\pa}{\pa y_j}}\,\om=0\,\,,\,\,\forall \,j\,\,.
\] 
Hence, we can write $\om=\sum_{j=1}^3A_j(x,y)\,dx_j$. On the other hand, the integrability condition $\om\wedge d\om=0$ implies that
\begin{equation}\label{eq:9}
A_i.\,\frac{\pa A_j}{\pa y_k}-A_j.\,\frac{\pa A_i}{\pa y_k}=0\,\,,\,\,\forall\,k=1,...,n-3\,\,,\,\,\forall\,1\le i<j\le 3\,\,.
\end{equation}
The relations in (\ref{eq:9}) imply that there exist $u\in\O^*_n$ and $B_1,B_2,B_3\in\O_3$ such that
\[
A_j(x,y)=u(x,y).\,B_j(x)\,\,,\,\,1\le j\le3\,\,\implies\,\,\om=u.\,\Pi^*(\wt\om)\,\,,
\]
where $\wt\om=\sum_{j=1}^3B_j(x)\,dx_j$ is integrable and $\wt\om\wedge (i_Z\,dx_1\wedge dx_2\wedge dx_3)=0$.
This contradicts proposition \ref{p:3}.
\qed

\begin{rem}\label{r:25}
{\rm Since the vector field $Z$ is polynomial, corollary \ref{c:2} can be globalised: in any rational manifold $M$ of dimension $n\ge3$ there are examples of codimension two foliations that are not contained in a codimension one foliation.}
\end{rem}

\section{Homogeneous foliations}\label{se:2}

\subsection{Homogeneous foliations}\label{se:21}
In this section we fix a coordinate system $(x_1,...,x_n)$ of $\C^n$. A $p$-form $\Om$ is said to be homogeneous of degree $m$ if its components are homogeneous polynomials of degree $m$. The radial vector field of $\C^n$ will be denoted by $R$:
\[
R=\sum_{j=1}^nx_j\,\frac{\pa}{\pa x_j}\,\,.
\] 
The $p$-form $\Om$ is said to be {\it dicritical} if $i_R\,\Om=0$. Otherwise we say that $\Om$ is non-dicritical.

If $\Om$ is a $p$-form of degree $m$ we have
\begin{equation}\label{eq:10}
L_R\,\Om=i_R\,d\Om+d\left(i_R\Om\right)=(m+p)\,\Om\,\,\,\,\,\text{(Euler's identity).}
\end{equation}

Next we state some useful results.

\begin{prop}\label{p:4}
Let $\eta$ be a homogeneous 2-form of degree $m$. Assume that $\eta$ is closed and $\eta^2=\eta\wedge\eta=0$. Then there exists a homogeneous integrable dicritical 1-form $\om$ such that $\eta=d\om$. In particular, $\om\wedge\eta=0$.
\end{prop}

\begin{rem}
{\rm The relation $\om\wedge\eta=0$ in the statement of proposition \ref{p:4} means that the leaves of the codimension two foliation $\fa_\eta$ are contained in the leaves of the codimension one foliation $\fa_\om$. In this case, we will say that {\it the foliation $\fa_\eta$ is contained in the foliation $\fa_\om$}.}
\end{rem}

{\it Proof of proposition \ref{p:4}.}
Set $\wt\om=i_R\eta$. Since $d\eta=0$, from Euler's identity (\ref{eq:10}) we get
\[
d\wt\om=d\left(i_R\eta\right)=L_R\eta=(m+2)\,\eta\,\,\implies\,\,\eta=d\om\,\,,\,\,\om=(m+2)^{-1}.\,\wt\om\,\,.
\]
Finally
\[
0=i_R(\eta\wedge\eta)=2\,i_R(\eta)\wedge\eta\,\,\implies\,\,\om\wedge d\om=(m+2)^{-1}.\,i_R(\eta)\wedge\eta=0\,\,\qed
\]

\begin{rem}
{\rm We would like to note that if the 1-form $\om$ is integrable, but not closed, then the 2-form $d\om$ is integrable and the foliation $\fa_{d\om}$ is contained in the foliation $\fa_\om$. For instance, if
\[
\om=f_1...,f_p\,\sum_j\la_j\,\frac{df_j}{f_j}\,\,,\,\,\la_j\in\C^*\,\,,\,\,\la_i\ne\la_j\,,\,\forall\,\,i\ne j\,\,,
\]
where $f_i$ is holomorphic $\forall i$, then
\[
d\om=d(f_1...f_p)\wedge\,\sum_j\la_j\,\frac{df_j}{f_j}\,\,.
\]

In this case, the leaves of $\fa_{d\om}$ are the connected components of the intersection of the levels of $f_1...f_p$ with the leaves of $\fa_\om$.

Let us mention that the foliation $\fa_{d\om}$ is contained in infinitely many codimension one foliations: all members of the family of foliations $\fa_{\om_\la}$, where
\[
\om_\la=f_1...,f_p\,\sum_j(\la_j+\la)\,\frac{df_j}{f_j}\,\,,\,\,\la\in\C\,\,.
\]}
\end{rem}

An analogous result in the homogeneous non-closed and non-dicritical case is the following:

\begin{prop}\label{p:5}
Let $\eta$ be a homogeneous integrable and non-dicritical 2-form on $\C^n$, $n\ge3$. Then the 1-form $\om:=i_R\,\eta$ is integrable. Moreover, the foliation $\fa_\eta$ is contained in the foliation $\fa_\om$.
\end{prop}

{\it Proof.}
We want to prove that $\om\wedge d\om=0$.
As in the proof of proposition \ref{p:4}, from $\eta^2=0$ we get
\[
0=i_R(\eta^2)=2\,i_R\eta\wedge\eta=2\,\om\wedge\eta\,\,\implies\,\,\om\wedge\eta=0\,\,.
\]
If $v$ is a constant vector field such that $f:=i_v\,\om\not\equiv0$ then
\[
\om\wedge\eta=0\,\,\implies\,\,i_v(\om\wedge\eta)=f.\,\eta-\om\wedge i_v\eta=0\,\,\implies\,\,\eta=\om\wedge\wt\om\,\,,
\]
where $\wt\om=f^{-1}.\,i_v\eta$. Therefore,
\[
\om=i_R\,\eta=-i_R\,\wt\om.\,\om\,\,\implies\,\,i_R\,\wt\om=-1\,\,.
\]
On the other hand, the integrability of $\eta$ implies that
\[
\om\wedge\wt\om\wedge d\om=0\,\,\implies\,\,0=i_R(\om\wedge\wt\om\wedge d\om)=-(i_R\wt\om).\,\om\wedge d\om+
\om\wedge\wt\om\wedge i_R\,d\om\,\,\implies
\]
\[
\om\wedge d\om=\wt\om\wedge\om\wedge i_R\,d\om\,\,.
\]
Finally, if $k=deg(\om)$ then (\ref{eq:10}) implies
\[
(k+1)\,\om=L_R\,\om=i_R\,d\om\,\,\implies\,\,i_R\,d\om\wedge\om=0\,\,\implies\,\,\om\wedge d\om=0\qed
\]
\vskip.1in
Let us state a result that will be important in what follows.
Let $\om$ be a homogeneous integrable 1-form on $\C^n$, $n\ge4$. We assume that $\om$ is dicritical of degree $k$ and $cod(Sing(\om))\ge2$.
Euler's identity (\ref{eq:10}) implies $i_R\,d\om=(k+1)\,\om$.
The form $\om$ induces a codimension one foliation $\fa$ on the space $\p^{n-1}$ whose singular set is the projectivisation of $Sing(\om)$;
$Sing(\fa)=\Pi(Sing(\om)\setminus\{0\})$, where $\Pi\colon\C^n\setminus\{0\}\to\p^{n-1}$ is the natural projection.

A {\it Kupka singularity} of $\fa$ is a point $p\in Sing(\fa)$ for which there is a local generator $\a$ of the germ $\fa_p$ with $d\a(p)\ne0$. The {\it Kupka set} of $\fa$ is, by definition, $K(\fa)=\{p\in Sing(\fa)\,|\,p$ is a Kupka singularity of $\fa\}$.
A {\it Kupka component} of $\fa$ is an irreducible component $K$ of $Sing(\fa)$ with $K\sub K(\fa)$.
It is known that a Kupka component of $\fa$ is a smooth sub-variety of codimension two along which the foliation is locally trivial
(see \S \ref{s:151}, [K] and [Me]).
If $\wt{K}$ is an irreducible component of $Sing(\om)$ such that $d\om(p)\ne0$ for all $p\in\wt{K}\setminus\{0\}$ then $\Pi(\wt{K}\setminus\{0\})$ is a Kupka component of $\fa$. The next statement resumes some results proven in [Ce-LN 1], [Br] and [CA-2]:

\begin{thma}\label{t:4}
Let $\fa$ and $\om$ be as above. If $\fa$ has a Kupka component $K$ then $K$ is a complete intersection of hypersurfaces, in homogeneous coordinates $(F=G=0)$. Moreover, if $degree(F)/degree(G)=p/q$, where $p,q\in\N$ and $(p,q)=1$, then $\fa$ is the algebraic pencil of hypersurfaces given by the rational function $\frac{F^q}{G^p}$, or equivalently by the 1-form $q\,F\,dG-p\,G\,dF$ (in homogeneous coordinates).
\end{thma}

\begin{rem}\label{r:23}
{\rm In fact, theorem \ref{t:4} says that we can choose $F$ and $G$ in such a way that $\om=q\,F\,dG-p\,G\,dF$. We would like to note also that the hypothesis $n\ge4$ is necessary; in dimension $n=3$ the statement is false.}
\end{rem}

\begin{Cor}\label{c:3}
Let $\om$ be an integrable homogenous and dicritical 1-form on $\C^n$, $n\ge4$. Then $0\in\C^n$ cannot be an isolated singularity of the 2-form $d\om$. 
\end{Cor}

{\it Proof.}
Suppose by contradiction that $Sing(d\om)=\{0\}$. In this case, all irreducible components of $\Pi(Sing(\om)\setminus\{0\})$ are contained in the Kupka set $K(\fa)$. Therefore, by theorem \ref{t:4} and remark \ref{r:23} we can suppose that $\om=q\,F\,dG-p\,G\,dF$, so that $d\om=(p+q)\,dF\wedge dG$ and $d\om$ is decomposable. However, this implies that $dim(Sing(d\om))\ge1$, a contradiction.
\qed

\subsection{Singularities of codimension two foliations}\label{se:22}
We would like to pose the following problem: 
\begin{prob}\label{pr:1}
Is there a germ of codimension two foliation with an isolated singularity at the origin of $\C^4$?
\end{prob}

First of all, in the case of dimension three there are such foliations. In fact, the codimension two foliations with an isolated singularity at the origin of $\C^3$ are generic.

Next, there are homogeneous codimension two distributions in $\C^4$ with an isolated singularity at the origin. An example, due to [K-N], is given by the decomposable but non-integrable 2-form 
\[
\te=x_3^2\,dx_2\wedge dx_3-x_1^2\,dx_3\wedge dx_1+(x_1\,x_2+x_3\,x_4)\,dx_1\wedge dx_2+
\]
\[
+[x_4^2\,dx_1+x_2^2\,dx_2+(x_1\,x_2-x_3\,x_4)\,dx_3]\wedge dx_4\,\,.
\]
The form $\te$ has an isolated singularity at $0\in\C^4$. It defines a distribution of 2-planes on $\C^4\setminus\{0\}$ because $\te\wedge\te=0$.

\vskip.1in

In fact, we don't know the answer of problem \ref{pr:1} in general, but the next statement contrasts with the previous example.
\begin{thm}\label{t:5}
Let $\eta$ be a homogeneous integrable 2-form on $\C^4$. Then $dim(Sing(\eta))\ge1$.
\end{thm}

{\it Proof.} We denote by $Z$ the rotational of $\eta$. We start by the case where $\eta$ is closed, which means $Z\equiv0$.
\vskip.1in
1- {\it $\eta$ is closed.} Let $\om=i_R\eta$. By proposition \ref{p:5}, $\om$ is integrable and by Euler's identity we have
\[
d\om=(m+2)\,\eta\,\,,
\] 
where $m$ is the degree of $\eta$. Since $dim(Sing(d\om))\ge1$ by corollary \ref{c:3}, we obtain the result in this case.
\vskip.1in
2- {\it $\eta$ is not closed, $Z\not\equiv0$.} Let us consider first the case where $i_R\eta=0$. In this case, since $0$ is an isolated singularity of $R$, by De Rham's division theorem there exists a homogeneous vector field $Y$ such that $\eta=i_R\,i_Y\,\nu$, $\nu=dx_1\wedge...\wedge dx_4$. Since $dim(Sing(R\wedge Y))\ge1$, we get $dim(Sing(\eta))\ge1$. In the same way, if $cod(Sing(Z))\ge3$, since $i_Z\,\eta=0$ then De Rham's division theorem implies that there exists a vector field $Y$ such that $\eta=i_Z\,i_Y\,\nu$ and again $dim(Sing(\eta))\ge1$.

Therefore, we can suppose that $cod(Sing(Z))\le2$ and $\om:=i_R\,\eta\not\equiv0$. In this case, by proposition \ref{p:5} the form $\om$ is integrable and induces a codimension one foliation on $\p^3$ of degree $\le m$, the degree of the coefficients of $\eta$. From Euler's identity (\ref{eq:10}) we get:
\begin{equation}\label{eq:111}
d\om=d\,i_R\,\eta=L_R\,\eta-i_R\,d\eta=(m+2)\,\eta-i_R\,i_Z\,\nu
\end{equation}
Since $i_Z\,\eta=0$ the above equality implies that $i_Z\,d\om=0$. Hence,
\[
L_Z\om=i_Z\,d\om+d\,i_Z\,\om=d\left(i_Z\,i_R\,\eta\right)=0\,\,\implies
\]
\begin{equation}\label{eq:11}
L_Z\,d\om=0
\end{equation}

Let us stablish a technical variant of corollary \ref{c:3}.

\begin{lemma}\label{l:21}
Let $\om$ be a homogeneous dicritical and integrable 1-form on $\C^n$, $n\ge4$. Suppose that $Sing(\om)$ contains an hypersurface $(h=0)$, $\om=h.\,\wt\om$, where $h(0)=0$ and $cod(Sing(\wt\om))\ge2$. Then $dim(Sing(d\om))\ge1$. 
\end{lemma}

{\it Proof.}
Suppose by contradiction that $Sing(d\om)=\{0\}$. Since
\[
d\om=dh\wedge \wt\om+h\,d\wt\om
\] 
we obtain that
\[
\wt\om(m)=0\,\,,\,\,m\ne0\,\,\implies\,\,d\om(m)\ne0\,\,\implies\,\,d\wt\om(m)\ne0\,\,.
\]
In particular, the singularities of $\wt\om$ in $\C^n\setminus\{0\}$ are of Kupka type and by theorem \ref{t:4} we have $\wt\om=k\,F\,dG-\ell\,G\,dF$, where $F$ and $G$ are homogeneous polynomials. Hence,
\[
d\om=dh\wedge\left(k\,F\,dG-\ell\,G\,dF\right)+(k+\ell)\,h\,dF\wedge dG
\]
and so $d\om$ vanishes on the set $(h=F=G=0)$, which in dimension $n\ge4$ has dimension $\ge1$.
\qed
\vskip.1in
Let us suppose now, by contradiction, that the 2-form $\eta$ has an isolated singularity at $0\in \C^4$. It follows from corollary \ref{c:3} and lemma \ref{l:21} that the 2-form $d\om$ vanishes at least on some straight line $L$ through $0\in\C^4$.
Moreover, from (\ref{eq:11}) the form $d\om$ is invariant by the local flow of $Z$. In particular, if we denote by $S$ the irreducible component of $Sing(d\om)$ that contains $L$, then $S$ is invariant by the local flows of $Z$ and of the radial vector field $R$. Hence $Z$ and $R$ are tangent to $S$. 
By Euler's identity (\ref{eq:111}), if $m\in S\setminus\{0\}$ then $d\om(m)=0$ and $\eta(m)\ne0$, imply that $R(m)$ and $Z(m)$ must be independent along $S$. Therefore, $dim(S)\ge2$. Since $Z\ne0$ along $S\setminus\{0\}$ ($S\cap Sing(Z)=\{0\}$) and we have supposed that $cod(Sing(Z))\le2$, of course we must have
\[
dim(S)=dim(Sing(Z))=2\,\,\implies\,\,S\setminus\{0\}\,\,\text{is smooth and connected}\,\,\implies
\]
its projectivisation in $\p^3$, $\Ga=\Pi(S\setminus\{0\})$, is a smooth curve. Let $\G$ be the one dimensional foliation on $\p^3$ defined in homogenous coordinates by the form $i_R\,i_Z\,\nu$. Since $R\wedge Z\ne0$ along $S\setminus\{0\}$ the curve $\Ga$ is an algebraic leaf of $\G$ such that $Sing(\G)\cap \Ga=\emp$. However, this is not possible by [Ln-So]: any algebraic curve invariant by a one dimensional foliation $\G$ of $\p^n$, $n\ge2$, must contain at least one singularity of $\G$ (see propositipon 2.4 in [LN-So]).
\qed

\begin{rem}
{\rm In general we don't know the answer of problem \ref{p:1}. However, a case in which $dim(Sing(\fa_\eta))\ge1$ is when there exists a germ of holomorphic vector field $Z$ such that $i_Z\eta=0$ and $cod(Sing(Z))\ge3$. Indeed, if this is true then by De Rham's division theorem we can write $\eta=i_Yi_Z\nu$. This implies that $dim(Sing(\fa_\eta))\ge1$, as in the argument of theorem \ref{t:5}. In particular, when $cod(Sing(rot(\eta))\ge3$ then $dim(Sing(\fa_\eta))\ge1$.}
\end{rem}

\subsection{Homogeneous integrable 2-forms of small degree}\label{se:23}

In this section $\eta$ will be a homogeneous integrable 2-form on $\C^n$, $n\ge4$.
Here we will describe with some detail the foliation $\fa_\eta$ when $0\le deg(\eta)\le 2$.

\vskip.1in

{\it The case $deg(\eta)=0$.} Here $\eta$ has constant coefficients, and so it is closed. Since $\eta^2=0$, by Darboux's theorem there exists a coordinate system $x=(x_1,...,x_n)$ such that $\eta=dx_1\wedge dx_2$. The leaves of $\fa_\eta$ are the level surfaces $(x_1=c_1,x_2=c_2)$.

\vskip.1in

{\it The case $deg(\eta)=1$.} We consider the 3-form with constant coefficients $d\eta$. We have two possilities according to $d\eta\not\equiv0$, or $d\eta\equiv0$.
If $d\eta\not\equiv0$ then it is integrable, by remark \ref{r:dom}. Hence, we can find coordinates $x=(x_1,...,x_n)$ such that $d\eta=dx_1\wedge dx_2\wedge dx_3$, so that,
\[
\eta=x_1\,dx_2\wedge dx_3+d\a\,\,,
\]
where $\a$ is a homogenous 1-form of degree two.
On the other hand, if $v_j=\frac{\pa}{\pa x_j}$ then $i_{v_j}\,d\eta=0$ for all $j\ge4$, then by remark \ref{r:dom} we get
\[
i_{v_j}\,\eta=0\,\,,\,\,\forall\,j\ge4\,\,\implies\,\,i_{v_j}d\a=0\,\,,\,\,\forall\,j\ge4\,\,\implies 
\]
\[
d\a=\sum_{1\le i<j\le\,3}A_{ij}(x)\,dx_i\wedge dx_j\,\,.
\]
Since $d\a$ is closed, we get
\[
\frac{\pa A_{ij}}{\pa x_k}=0\,\,,\,\,\forall\,j\ge4\,\,\implies\,\,A_{ij}(x)\,\,\text{depends only of $x_1$, $x_2$ and $x_3$.}
\]

In particular, we can assume that $\a$ does not depend of the variables $x_4,...,x_n$:
\[
\a=\sum_{j=1}^3A_j(x_1,x_2,x_3)\,dx_j\,\,,
\]
where $A_1$, $A_2$ and $A_3$ are homogeneous polynomials of degree two.
The foliation $\fa_\eta$ is therefore the pull-back of a foliation on $\C^3$ defined by a linear vector field
\[
L=\left(x_1+A_{23}(x_1,x_2,x_3)\right)\frac{\pa}{\pa x_1}-A_{13}(x_1,x_2,x_3)\frac{\pa}{\pa x_2}+A_{12}(x_1,x_2,x_3)\frac{\pa}{\pa x_3}\,\,,
\]
where $A_{ij}=\frac{\pa A_j}{\pa x_i}-\frac{\pa A_i}{\pa x_j}$.

If $d\eta\equiv0$ then by proposition \ref{p:4} we have $\eta=d\om$, where $\om=\frac{1}{3}\,i_R\eta$ is dicritical and integrable (proposition \ref{p:5}). The 1-form $\om$ defines in homogeneous coordinates a foliation of degree one on $\p^{n-1}$. 
It is known in this case that $\om$ has an integrating factor: there exists a homogeneous polynomial $P$ of degree three such that $\frac{\om}{P}$ is closed (cf. [J] or [Ce-LN]).
As a consequence,
\[
\eta=d\om=\frac{dP}{P}\wedge\om\,\,,
\]
and so $\fa_\eta$ is the intersection of two foliations, $\fa_\om$ and $\fa_{dP}$.
In this situation there are two generic cases (cf. [Ce-LN]):

\vskip.1in

{\it $1^{st}$ case.} $P=L.Q$, where $L$ is linear and $Q$ generic of degree two. Here we have
\[
\om=\frac{1}{3}\left(L\,dQ-2\,Q\,dL\right)\,\,\implies\,\,\eta=dL\wedge dQ\,\,.
\]  
The polynomial map  $(L,Q)\colon\C^n\to\C^2$ is a first integral of $\fa_\eta$. 
\vskip.1in
{\it $2^{nd}$ case.} $P=L_1.\,L_2.\,L_3$, where $L_1$, $L_2$ and $L_3$ are linear and independent. We can assume that $P=x_1.\,x_2.\,x_3$. Here we have
\[
\om=x_1.\,x_2.\,x_3\left(\la_1\,\frac{dx_1}{x_1}+\la_2\,\frac{dx_2}{x_2}+\la_3\,\frac{dx_3}{x_3}\right)\,\,,\,\la_1+\la_2+\la_3=0\,\,.
\]
In this case,
\[
\eta=d(x_1.\,x_2.\,x_3)\wedge\left(\la_1\,\frac{dx_1}{x_1}+\la_2\,\frac{dx_2}{x_2}+\la_3\,\frac{dx_3}{x_3}\right)
\]
and $\fa_\eta$ has two first integrals, the first one algebraic $x_1.\,x_2.\,x_3$ and the second one in general "liouvillian" $x_1^{\la_1}.\,x_2^{\la_2}.\,x_3^{\la_3}$, but maybe algebraic when $[\la_1:\la_2:\la_3]\in\p^2_\Z$.

\vskip.1in
{\it The case $deg(\eta)=2$.} This case is more difficult. We again distinguish the two cases, $d\eta\equiv0$ and $d\eta\not\equiv0$.
If $d\eta\equiv0$ then by proposition \ref{p:4} there exists a dicritical homogeneous 1-form of degree three $\om$ such that $\eta=d\om$.
The form $\om$ is the homogeneos expression of some codimension one foliation on $\p^{n-1}$ of degree two.
According to [Ce-LN] the space of such foliations has six irreducible components:
\vskip.1in
$1^{st}$: $\ov{R(2,2)}$. Here the generic member has a rational first integral of the form $\frac{P}{Q}$, where $P$ and $Q$ are quadrics. In this case, $\om=\frac{1}{2}(P\,dQ-Q\,dP)$ and $\eta=dP\wedge dQ$. The foliation $\fa_\eta$ has the first integral $(P,Q)\colon \C^n\to\C^2$.
\vskip.1in
$2^{nd}$: $\ov{R(1,3)}$. Here the generic member has a rational first integral of the form $\frac{C}{L^3}$, where $C$ is a cubic and $L$ linear. In this case, $\om=\frac{1}{4}(L\,dC-3\,C\,dL)$ and $\eta=dL\wedge dC$. The foliation $\fa_\eta$ has the first integral $(L,C)\colon \C^n\to\C^2$.
\vskip.1in
$3^{rd}$: $\ov{L(1,1,1,1)}$. Here the generic member can be expressed in homogeneous coordinates by a the 1-form
\[
\om=L_1.\,L_2.\,L_3.\,L_4\left(\la_1\,\frac{dL_1}{L_1}+\la_2\,\frac{dL_2}{L_2}+\la_3\,\frac{dL_3}{L_3}+\la_4\,\frac{dL_4}{L_4}\right)\,\,,
\]
where the $L_{j}$ is linear, $\la_j\in\C^*$, $1\le j\le4$, and $\sum_j\la_j=0$. In this case, we have
\[
\frac{1}{L_1.\,L_2.\,L_3.\,L_4}\,\eta=\sum_i\frac{dL_i}{L_i}\,\wedge\,\sum_j\la_j\frac{dL_j}{L_j}
\]
and $\fa_\eta$ is the intersection of the two foliations $\fa_{d(L_1L_2L_3L_4)}$ and $\fa_\om$.
\vskip.1in
$4^{th}$: $\ov{L(1,1,2)}$. Here the generic member can be expressed as
\[
\om=L_1.\,L_2.\,Q\left(\la_1\,\frac{dL_1}{L_1}+\la_2\,\frac{dL2}{L_2}+\la\,\frac{dQ}{Q}\right)\,\,,
\]
where $L_1$ and $L_2$ are linear, $Q$ a quadric and $\la_1+\la_2+2\la=0$. Again the foliation $\fa_\eta$ is the intersection of two others: the foliations $\fa_{d(L_1L_2\,Q)}$ and $\fa_\om$.
\vskip.1in
$5^{th}$: $\ov{E(n-1)}$. Here the generic member has a first integral of the form $F=\frac{C^2}{Q^3}$ where $C$ is a cubic and $Q$ a quadric. The form $C.\,Q\,\frac{dF}{F}=2\,Q\,dC-3\,C\,dQ$ has a linear factor. It is proved in [Ce-LN] that in some homogeneous coordinate system $x=(x_1,...,x_4,...,x_n)$ we can write
\[
C=x_3\,x_4^2-x_1\,x_2\,x_4+\frac{x^3_1}{3}\,\,\,\,\text{and}\,\,\,\,Q=x_2\,x_4-\frac{x_1^2}{2}
\]
and the linear factor is $x_4$. In these coordinates we have $\om=\frac{1}{x_4}(2\,Q\,dC-3\,C\,dQ)$ and
\[
\eta=d\left(\frac{C\,Q}{x_4}\right)\wedge \left(2\frac{dC}{C}-3\frac{dQ}{Q}\right)\,\,\implies
\]
the foliation $\fa_\eta$ is the intersection of $\fa_{d(CQ/x_4)}$ and $\fa_\om$.
\vskip.1in
$6^{th}$: $S(2,n)$. Here the foliation induced by $\om$ in $\p^{n-1}$ is a linear pull-back of a degree two foliation on $\p^2$. This means that there exist homogeneous coordinates $x=(x_1,x_2,x_3,...,x_n)$ on $\C^n$ and homogeneous polynomials of degree three $P$, $Q$ and $R$, depending only of $x_1,x_2,x_3$, such that $x_1\,P+x_2\,Q+x_3\,R\equiv0$ and
\[
\om=P(x_1,x_2,x_3)\,dx_1+Q(x_1,x_2,x_3)\,dx_2+R(x_1,x_2,x_3)\,dx_3\,\,\implies
\]
\[
\eta=dP\wedge dx_1+dQ\wedge dx_2+dR\wedge dx_3\,\,.
\]
In other words, $\fa_\eta$ is the pull-back by a projection $x\in\C^n\mapsto(x_1,x_2,x_3)\in\C^3$, of a homogeneous foliation of degree two and codimension two in $\C^3$.

Let us mention that in the above case, all leaves of $\fa_\eta$ are rulled: they contain the fibers of the projection $(x_1,...,x_n)\mapsto (x_1,x_2,x_3)$. On the other hand, in general a foliation of degree two on $\p^2$ has no algebraic leaves. Therefore, in general the leaves of the foliations $\fa_\om$ and $\fa_\eta$ are Zariski dense.
\vskip.1in
Let us suppose now $d\eta\not\equiv0$.
We first consider the case where $n=4$ and the rotational of $\eta$, $X:=rot(\eta)$, satisfies $cod(Sing(X))\ge3$.
\vskip.1in
\begin{lemma}\label{l:22}
In the above situation we have
\[
\eta=i_Y\,i_X\nu\,\,,\,\,\nu=dz_1\wedge dz_2\wedge dz_3\wedge dz_4\,\,,
\]
where $Y$ is a linear vector field satisfying $[Y,X]=\la X$, $\la=1-tr(Y)$. Moreover,
\begin{itemize}
\item[(a).] If $X$ is not nilpotent then $\la=0$ and $tr(Y)=1$. In particular, $X$ and $Y$ commute.
\item[(b).] If $X$ is nilpotent and $\la\ne0$ then after a linear change of variables we have
\[
X=z_1\frac{\pa}{\pa_{z_2}}+z_2\frac{\pa}{\pa_{z_3}}+z_3\frac{\pa}{\pa_{z_4}}
\]
and
\[
Y=\rho z_1\frac{\pa}{\pa z_1}+(\rho-\la)z_2\frac{\pa}{\pa z_2}+(\rho-2\la)z_3\frac{\pa}{\pa z_3}+(\rho-3\la)z_4\frac{\pa}{\pa z_4}
\]
where $4\rho-5\la=1$.
\end{itemize}
\end{lemma}

{\it Proof.}
Recall that $d\eta=i_X\nu$, $\nu=dx_1\wedge...\wedge dx_4$. Since $\eta$ is homogeneous of degree two, $X$ is a linear vector field with $tr(X)=0$.
Since $cod(sing(X))\ge3$, by the division theorem, there exists another linear vector field $Y$ such that
\[
\eta=i_Y\,i_X\,\nu\,\,\implies\,\,L_Y\,\eta=i_Y\,d\eta=\eta\,\,\implies\,\,L_Y\,d\eta=d\eta\,\,.
\]
The last relation implies
\[
d\eta=i_X\,\nu=L_Y(i_X\,\nu)=i_{[Y,X]}\,\nu+i_X\,L_Y\,\nu=i_{[Y,X]}\,\nu+tr(Y)\,i_X\,\nu\,\,\implies
\]
\[
[Y,X]=(1-tr(Y))\,X:=\la\,X\,\,,\,\,\la=1-tr(Y)\,\,.
\]

Consider $X$ as a derivation $X=\sum_{j=1}^4X_j\,\frac{\pa}{\pa x_j}$. Since $X$ is linear the $k^{th}$ power operator $X^k$, $k\ge2$, is also a derivation $X^k=\sum_{j=1}^4X_j^k\,\frac{\pa}{\pa x_j}$. Moreover, if the eigenvalues of $X$ are $\la_1,...,\la_4$ then the eigenvalues of $X^k$ are $\la_1^k,...,\la_4^k$.

As a derivation, the relation $[Y,X]=\la\,X$ can be written as $Y.\,X-X.\,Y=\la\,X$.
It implies that $Y.\,X^k-X^k.\,Y=k\,\la\,X^k$, for all $k\ge1$.
The proof is by induction on $k\ge1$. Let us assume, by induction, that $Y.\,X^k-X^k.\,Y=k\,\la\,X^k$ for some $k\ge1$.
Then
\[
\left.
\begin{matrix}
Y.\,X^k-X^k.\,Y=k\,\la\,X^k\\
Y.\,X-X.\,Y=\la\,X\\
\end{matrix}
\right\}\,\,\implies\,\,
\left.
\begin{matrix}
Y.\,X^{k+1}-X^k.\,Y.\,X=k\,\la\,X^{k+1}\\
X^k.\,Y.\,X-X^{k+1}.\,Y=\la\,X^{k+1}\\
\end{matrix}
\right\}\,\,\implies
\]
\[
Y.\,X^{k+1}-X^{k+1}.\,Y=(k+1)\,\la\,X^{k+1}\,\,.
\]
Therefore, $Y.\,X^k-X^k.\,Y=k\,\la\,X^k$ for all $k\ge1$. If $\la\ne0$ then $tr(X^k)=k^{-1}.\,\la^{-1}.\,tr(Y.\,X^k-X^k.\,Y)=0$ for all $k\ge1$. However, this implies that $\la_1=...=\la_4=0$ and that $X$ is nilpotent.
Therefore, if $X$ is not nilpotent then $\la=0$, which proves (a).

Assume that $X$ is nilpotent and $\la\ne0$. After a linear change of variables we can assume that
$X=z_1\frac{\pa}{\pa z_2}+z_2\frac{\pa}{\pa z_3}+z_3\frac{\pa}{\pa z_4}$. As a derivation we have $[Y,X]=Y.\,X-X.\,Y=\la\,X$. If we apply both members in $z_1$ then we get $X(Y(z_1))=0$ which implies that $Y(z_1)$ is an eigenvector of $X$: $Y(z_1)=\rho\,z_1$, $\rho\in\C$.
When we apply in $z_2$ then $(Y.X-X.Y)(z_2)=Y(z_1)-X(Y(z_2))=\la\,z_1\,\,\implies$
\[
X(Y(z_2))=(\rho-\la)\,z_1\,\,\implies\,\,Y(z_2)=(\rho-\la)\,z_2+a\,z_1\,\,,\,\,a\in\C\,\,.
\]
By a similar argument we obtain
\[
Y(z_3)=(\rho-2\,\la) z_3+a\,z_2+b\,z_1\,\,\text{and}\,\,Y(z_4)=(\rho-3\,\la)z_4+a\,z_3+b\,z_2+c\,z_1\,\,,\,\,b\,,\,c\in\C\,\,.
\]
In particular, the eigenvalues of $Y$ are $\rho$, $\rho-\la$, $\rho-2\,\la$ and $\rho-3\,\la$. Hence, $tr(Y)=4\,\rho-6\,\la$ and since $\la=1-tr(Y)$ we get the relation $4\rho-5\la=1$.
The eigenvalues of $Y$ are two by two distinct and so it is diagonalizable. Therefore, $Y$ has an eigenvector $w=z_4+\a\,z_3+\be\,z_2+\g\,z_1$, where $Y(w)=(\rho-3\la)w$. Set $z:=X(w)=z_3+\a\,z_2+\be\,z_1$, $y:=X(z)=z_2+\a\,z_1$ and $x:=X(y)=z_1$. Finally,
\[
Y(X(w))-X(Y(w))=\la\,X(w)\,\,\implies\,\,Y(z)=(\rho-2\la)\,z\,\,.
\]
Similarly $Y(y)=(\rho-\la)\,y$ and $Y(x)=\rho\,x$. This finishes the proof.
\qed
\vskip.1in

In case (a) of lemma \ref{l:22}, where $[X,Y]=0$, the vector fields $X$ and $Y$ generate an action of $\C^2$ on $\C^4$. We will assume the generic case, in which the $X$ and $Y$ are diagonalizable in the same basis of $\C^4$. This means that after a linear change of variables we can assume that $X=\sum_{j=1}^4\la_j\,z_j\frac{\pa}{\pa z_j}$ and $Y=\sum_{j=1}^4\mu_j\,z_j\frac{\pa}{\pa z_j}$, where $\sum_j\la_j=0$ and $\sum_j\mu_j=1$.
We will assume also that $\la_i\,\mu_j-\la_j\,\mu_i\ne0$ if $i\ne j$.
In this case, we have
\begin{equation}\label{eq:13'}
\eta=i_Yi_X\nu=z_1\,z_2\,z_3\,z_4\sum_{i<j}\rho_{ij}\frac{dz_i\wedge dz_j}{z_i\,z_j}\,,
\end{equation}
where $\nu=dz_1\wedge dz_2\wedge dz_3\wedge dz_4$, $\rho_{ij}=\pm(\la_k\mu_\ell-\la_\ell\mu_k)$ and $\{k,\ell\}=\{1,2,3,4\}\setminus\{i,j\}$. In particular, $f:=z_1z_2z_3z_4$ is an integrating factor of $\eta$: $d\left(\frac{\eta}{f}\right)=0$.

\vskip.1in

As we will see next, in case (b) of lemma \ref{l:22} the form $\eta$ also has an integrating factor.
Inded, in case (b) we can write $Y=\rho\,R-\la\,S$, where $R=\sum_{j=1}^4z_j\frac{\pa}{\pa z_j}$ is the radial vector field and $S=\sum_{j=2}^4(j-1)z_j\,\frac{\pa}{\pa z_j}$. 
In particular, $S,R,X$ generate a Lie algebra of linear vector fields with the relations $[S,R]=[R,X]=0$ and $[S,X]=-X$.
Set $\a:=i_S\,i_X\nu$ and $\be:=i_R\,i_X\nu$, so that $\eta=\rho.\,\be-\la.\,\a$. Let $f:=g.h/z_1$, where
\begin{equation}\label{eq:14'}
g=z_2^3-3\,z_1\,z_2\,z_3+3\,z_1^2\,z_4\,\,\text{and}\,\,h=z_2^2-2\,z_1\,z_3\,\,.
\end{equation}
The reader can check directly that $df\wedge\a=f.\,d\a$ and $df\wedge \be=f.\,d\be$, which is equivalent to $d\left(\frac{1}{f}\,\a\right)=0$ and 
$d\left(\frac{1}{f}\,\be\right)=0$ and this implies $d\left(\frac{1}{f}\,\eta\right)=0$.
It follows also that $\frac{1}{f}\,\a$ and $\frac{1}{f}\,\be$ are logarithmic 2-forms with pole divisor $z_1.\,g.\,h$; they belong to the vector space generated by $\frac{dg\wedge dh}{gh}$, $\frac{dh\wedge dz_1}{hz_1}$ and $\frac{dz_1\wedge dg}{z_1g}$. 
In fact, the reader can check directly that
\[
\frac{z_1\a}{gh}=\left(\frac{1}{3}\frac{dg}{g}-\frac{1}{2}\frac{dh}{h}\right)\wedge \frac{dz_1}{z_1}
\,\,\,\text{and}\,\,\,
\frac{z_1\be}{gh}=\frac{1}{6}\frac{dh\wedge dg}{hg}+\frac{z_1\a}{gh}
\]
so that
\begin{equation}\label{eq:13}
\frac{z_1\eta}{gh}=A\,\frac{dh\wedge dg}{hg}+B\,\frac{dg\wedge dz_1}{gz_1}+C\,\frac{dz_1\wedge dh}{z_1h}\,\,,
\end{equation}
where $A=\rho/6$, $B=(\rho-\la)/3$ and $C=(\rho-\la)/2$.
\begin{rem}
{\rm As we have seen in proposition \ref{p:5} the form $\om:=i_R\eta$ is integrable; $\om\wedge d\om=0$. Since $i_Rd\om=4\,\om$, if $\om\not\equiv0$ then $d\om\not\equiv0$. We would like to observe that for every $s\in\C$ the form $\eta_s:=\eta+s\,d\om$ is integrable. Let us prove this fact.

First of all Euler's identity implies that
\begin{equation}\label{eq:14}
d\om=d\,i_R\eta=4\,\eta-i_R\,d\eta=4\,\eta-i_Ri_X\nu\,\,\implies
\end{equation}
\[
d\om\wedge \eta=-i_Ri_X\nu\wedge i_Yi_X\nu=0\,\,\implies
\]
\[
\eta_s^2=(\eta+s\,d\om)^2=0\,\,,
\]
because $\eta^2=d\om^2=0$. On the other hand $rot(\eta_s)=rot(\eta)=X$ and from (\ref{eq:14}) we get $i_X\eta_s=0$.}
\end{rem}

Let us consider now the case $\C^n$, $n\ge5$.

\begin{prop}\label{p:6}
Let $\eta$ be an integrable homogeneous of degree two form on $\C^n$, $n\ge5$. Assume that there exists a 4-plane $\Si_o$ with $0\in\Si_o\sub\C^n$ such that
\begin{itemize}
\item[(a).] The singular set of $rot(\eta|_{\Si_o})$ has codimension $\ge3$.
\item[(b).] $\eta|_{\Si_o}$ is generic, as in {\rm(\ref{eq:13'})} or {\rm (\ref{eq:13})}.
\end{itemize}
Then there exists a linear coordinate system $(z_1,z_2,z_3,z_4,...,z_n)\in\C^n$ where $\eta$ can be written, either as in {\rm(\ref{eq:13'})}, or as in {\rm(\ref{eq:13})}. In particular, $\eta$ depends only of four variables. 
\end{prop}

{\it Proof.}
We will assume first that $i_R\,\eta\not\equiv0$, $R=\sum_{j=1}^nz_j\,\frac{\pa}{\pa z_j}$. 

Consider the integrable 1-form $\om=i_R\eta$. 
Then $\om$ represents in homogeneous coordinates a foliation of degree two on $\p^{n-1}$, denoted by $\fa_\om$. 
Let $\Si_o$ be a 4-plane as in the hypothesis.
We can assume that $\Si_o=\{(x,0)\,|\,x=(z_1,z_2,z_3,z_4)\}=\C^4$ and that $\eta|_{\Si_o}$ is either as in (\ref{eq:13'}), or as in (\ref{eq:13}). 
Since $R$ is tangent to $\Si_o$ we obtain $\om|_{\Si_o}=i_{R_4}(\eta|_{\Si_o})$, where $R_4=\sum_{j=1}^4z_j\,\frac{\pa}{\pa z_j}$.
Now we observe the following:
\begin{itemize}
\item{1$^{st}$.} If $\eta|_{\Si_o}$ is as in (\ref{eq:13'}) then $\fa_\om|_{\Si_o}=\fa_{(\om|_{\Si_o})}\in L(1,1,1,1)$.
\item{2$^{nd}$.} If $\eta|_{\Si_o}$ is as in (\ref{eq:13}) then $\fa_\om|_{\Si_o}=\fa_{(\om|_{\Si_o})}\in E(3)$.
\end{itemize}
It was proved in [Ce-Ln] that in the 1$^{st}$ case then $\fa_\om\in L(1,1,1,1)$ and $\om$ has an integrating factor of the form $F=\ell_1.\ell_2.\ell_3.\ell_4$, where $\ell_j$ is linear $1\le j\le4$. In some linear coordinate system $z=(z_1,...,z_n)$ we can assume that $F=z_1.z_2.z_3z_4$ and
\[
\om=z_1.z_2.z_3.z_4\sum_{j=1}^4\la_j\frac{dz_j}{z_j}\,\,\text{, where}\,\,\sum_j\la_j=0\,\,.
\]
On the other hand, in the 2$^{nd}$ then $\fa_\om\in E(n-1)$ and in some linear coordinate system we have $z_1.\,\om=3g.\,dh-2h.\,dg$, where $g=g(z_1,z_2,z_3,z_4)$ and $h=h(z_1,z_2,z_3)$ are as in (\ref{eq:14'}). 

In both cases we can write $\om=\sum_{j=1}^4A_j(z_1,z_2,z_3,z_4)\,dz_j$. In other words, $\om$ depends only of $z_1,...,z_4$ and $dz_1,...,dz_4$.
We assert that the same is true for $\eta$:
\[
\eta=\sum_{1\le i<j\le 4}A_{ij}(z_1,z_2,z_3,z_4)\,dz_i\wedge dz_j\,\,.
\]
Let us prove this fact. First of all, from $\eta^2=0$, we get
\begin{equation}\label{eq:17}
\om\wedge\eta=\frac{1}{2}i_R(\eta\wedge\eta)=0
\end{equation}
Let $e_j=\frac{\pa}{\pa z_j}$, $j=1,...,n$. Relation (\ref{eq:17}) and $i_{e_j}\om=0$, $j\ge5$, imply that $\om\wedge i_{e_j}\eta=0$, $\forall\,\,j\ge5$. Since $cod(Sing(\om))=2$ in all cases, from the division theorem we obtain $i_{e_j}\eta=f_j.\,\om$, where $f_j$ is holomorphic. However, since $\om$ is homogeneous of degree three and $\eta$ homogeneous of degree two, we must have $i_{e_j}\eta=0$ for all $j\ge5$.
Therefore, we can write $\eta=\sum_{1\le i<j\le4}B_{ij}(z)\,dz_i\wedge dz_j$, where $B_{ij}$ is homogeneous of degree two, $1\le i<j\le4$.
We assert that $B_{ij}$ depends only of $z_1,z_2,z_3,z_4$, or equivalently that
\[
L_{e_k}\eta=\sum_{1\le i<j\le4}\frac{\pa B_{ij}}{\pa z_k}\,dz_i\wedge dz_j=0\,\,,\,\,\forall\,\,k\ge5\,\,.
\]
First of all, since $\om=i_R\eta$ and $\om$ and $\eta$ depend only of $dz_1,dz_2,dz_3,dz_4$, we can write
\[
\om=i_{R_4}\eta\,\,,\,\,R_4=\sum_{j=1}^4z_j\,\frac{\pa}{\pa z_j}\,\,.
\]
Since $[R_4,e_k]=0$ for $k\ge5$, we get from the relation $\om=i_{R_4}\eta$ that $i_{R_4}\,L_{e_k}\eta=0$ and $i_{R_4}\,L_{e_k}L_{e_\ell}\eta=0$, $\forall k,\ell\ge5$. Since $L_{e_k}L_{e_\ell}\eta$ has constant coefficients, this last relation implies that $L_{e_k}L_{e_\ell}\eta=0$, $\forall k,\ell\ge5$. In particular, we can write $B_{ij}(z)=A_{ij}(z_1,...,z_4)+\sum_{k\ge5}z_k\,C_{ijk}(z_1,...,z_4)$, $\forall 1\le i<j\le4$, or equivalently $\eta=\eta_o+\sum_{k\ge5}z_k\,\eta_k$, where
\[
\eta_o=\sum_{1\le i<j\le4}A_{ij}(z_1,...,z_4)dz_i\wedge dz_j\,\,\text{and}\,\,\eta_k=L_{e_k}\eta=\sum_{1\le i<j\le4}C_{ijk}(z_1,...,z_4)dz_i\wedge dz_j\,\,.
\]
Now, from the relation $\eta^2=0$ we obtain
\[
\eta_o^2+2\,\sum_{k\ge5}z_k\eta_o\wedge\eta_k+\sum_{k,\ell\ge5}z_k\,z_\ell\,\eta_k\wedge\eta_\ell\equiv=0\,\,\implies
\]
$\eta_o^2=0$ and
\begin{equation}\label{eq:18}
\eta_o\wedge\eta_k=0\,\,,\,\,\forall\,k\ge5\,\,.
\end{equation}

Recall that we have assumed that $\eta|_{\Si_o}$ is like in (\ref{eq:13'}) or in (\ref{eq:13}), but since $\Si_o=\{(z_1,...,z_4,0...,0)\,|\,(z_1,...,z_4)\in\C^4\}$, we have $\eta|_{\Si_o}=\eta_o|_{\Si_o}=\eta_o$.
In other words, we can consider $\eta_o$ and $\eta_{k}$, $k\ge5$, as 2-forms on the 4-plane $\Si_o$.
However, as we have seen, there is a linear vector field $X=rot(\eta_o)$ on the plane $\Si_o$ such that $cod(Sing(X))\ge3$ and $i_X\eta_o=0$.
Recall also that $i_{R_4}\eta_k=0$. Hence, 
from (\ref{eq:18}) we obtain
\[
\eta_o\wedge i_X\eta_k=0\,\,\implies\,\,\om\wedge i_X\eta_k=i_{R_4}\eta_o\wedge i_X\eta_k=0\,\,,\,\,\forall\,k\ge5\,\,.
\]
As before, this implies $i_X\eta_k=0$. But, since $\eta_k$ is homogeneous of degree one this together $i_{R_4}\eta_k=0$ is possible only if $\eta_k=0$. 

\vskip.1in

It remains to study the case $i_R\eta=0$. Observe first that, since $i_R\eta=0$, we can consider $\eta$ as the homogeneous expression of a foliation of degree one on $\p^{n-1}$.
Such foliations were classified in [L-P~T] in any codimension. There are two irreducible components in general. In the case of codimansion two foliations, in one of the components the typical foliation is a linear pull-back of a 1-dimensional foliation on $\p^3$ of degree one.
In the other component, the typical foliation is the pull-back of a one dimensional foliation on the weighted projective space $\p^3_{(2,1,1)}$.

In our case we are assuming that $\eta|_{\Si_o}$ is like in (\ref{eq:13'}) or (\ref{eq:13}).  
We assert that $\eta|_{\Si_o}$ cannot be as in (\ref{eq:13}).

Indeed, choose coordinates $z=(u,v)\in\C^4\times\C^{n-4}$, $u=(z_1,...,z_4)$ and $v=(z_5,...,z_n)$, such that $\Si_o=(v=0)$.
Note that $i_{R_4}\eta|_{\Si_o}=0$, where $R_4=\sum_{j=1}^4z_j\frac{\pa}{\pa z_j}$. On the other hand, in case (\ref{eq:13}) we have $\eta|_{\Si_o}=\rho\,i_{R_4}i_{X_o}\nu-\la\,i_Si_{X_o}\nu$,
where $X_o=rot(\eta|_{\Si_o})=z_1\frac{\pa}{\pa z_2}+z_2\frac{\pa}{\pa z_3}+z_3\frac{\pa}{\pa z_4}$ and $S=\sum_j(j-1)z_j\frac{\pa}{\pa z_j}$.
Therefore,
\[
i_{R_4}\eta|_{\Si_o}=-\la\,i_{R_4}i_Si_{X_o}\nu\not\equiv0\,,
\]
because $\la\ne0$ and $R_4\wedge S\wedge X_o\ne0$, as the reader can check.

Finally, if $\eta|_{\Si_o}$ is like in (\ref{eq:13'}) then it can be considered as the homogeneous expression of a 1-dimensional foliation of degree one on $\p^3$. Therefore, by the classificationin [L-P-T], the form $\eta$ is the linear pull-back of $\eta|_{\Si_o}$ by a linear map $T\colon\C^n\to\C^4$. This finshes the proof of the proposition.
\qed 

\vskip.1in

Let us complete the study of in $\C^4$. It remains to consider the cases in which the singular set of $X=rot(\eta)$ has codimension one or two.

\vskip.1in

{\it 1$^{st}$ case}: $\eta$ is dicritical, $i_R\eta=0$. We assert that in this case we have: $\eta=\frac{1}{4}\,i_R\,i_X\,\nu$, $\nu=dz_1\wedge...\wedge dz_4$. In particular, $\fa_\eta$ is defined by a commutative action.

In fact, $i_R\eta=0$ implies that $\eta=i_R\mu$, where $\mu$ is a 3-form homogeneous of degree one.
In particular, there exits a linear vector field $Y$ such that $\mu=i_Y\nu$, so that $\eta=i_Ri_Y\nu$.
On the other hand, Euler's identity implies
\[
i_X\,\nu=d\eta=d\,i_R\mu=L_R\mu-i_R\,d\mu=4\,\mu-i_R\,d\mu=4\,i_Y\,\nu-i_R\,d\mu\,\,.
\]  
Since $d\mu$ is homogeneous of degree zero, then $d\mu=\rho\,\nu$, where $\rho\in\C$. Therefore, the above relation implies that
\[
X=4\,Y-\rho\,R\,\,\implies\,\,\eta=\frac{1}{4}\,i_R\,i_X\,\nu
\]
as asserted.

\vskip.1in

{\it 2$^{nd}$ case}: $\eta$ is non-dicritical, $\om:=i_R\eta\not\equiv0$. We have two sub-cases:

\vskip.1in

{\it Case 2.1}: $Sing(X)$ has codimension one. In this case, we can write $X=H.\,Y$, where $H$ is linear and $Y$ a constant vector field.
After a linear change of variables we can assume that $Y=\frac{\pa}{\pa z_4}$, and so
\[
d\eta=H\,dz_1\wedge dz_2\wedge dz_3\,\,\implies\,\,\frac{\pa H}{\pa z_4}=0\,\,\implies\,\,H=H(z_1,z_2,z_3)\,\,.
\]
Since $i_X\eta=0$ we can write
\[
\eta=A_1\,dz_2\wedge dz_3+A_2\,dz_3\wedge dz_1+A_3\,dz_1\wedge dz_2=i_Z\,dz_1\wedge dz_2\wedge dz_3\,\,, 
\]
where $Z=\sum_{j=1}^3A_j\frac{\pa}{\pa z_j}$. From $d\eta=H\,dz_1\wedge dz_2\wedge dz_3$, we get $\frac{\pa A_j}{\pa z_4}=0$, $1\le j\le 4$, and $\sum_{j=1}^3\frac{\pa A_j}{\pa z_j}=H$. Therefore, $\fa_\eta$ is a linear pull-back of a degree two homogeneous one dimensional foliation on $\C^3$: the foliation defined by $Z$ on $\C^3$.

\begin{rem}
{\rm In general $\eta$ has no rational integrating factor in this case. 
Indeed, if it had a rational integrating factor, say $f=f_1/f_2$, $d\frac{1}{f}\eta=0$, then $df\wedge\eta=f.\,d\eta$, which implies that $\frac{\pa f}{\pa z_4}=0$ and $Z(f_j)=g_j.\,f_j$, $j=1,2$. In other words, the foliation defined by $Z$ has least one invariant homogeneous hypersurface. However, this is not true in general. In fact, consider the 1-form
\[
\om=i_Ri_Z\,dz_1\wedge dz_2\wedge dz_3\,\,.
\]
The form $\om$ can be considered as the homogeneous expression of a degree two foliation on $\p^2$, say $\G$. Also, $\om$ has the same invariant homogeneous hypersurfaces as $\eta$. A homogeneous invariant hypersurface for $\eta$ gives origin to an algebraic invariant curve for $\G$. However, it is known that a generic foliation of degree two on $\p^2$ has no invariant algebraic curve (cf. [J] and [LN 1]).}
\end{rem}

\vskip,1in

{\it Case 2.2}: $Sing(X)$ has codimension two. Since $X$ is linear and $tr(X)=0$, it corresponds to a rank two $4\times4$ matrix with vanishing trace. There are three possible Jordan canonical forms: $z_1\,\frac{\pa}{\pa z_1}-z_2\,\frac{\pa}{\pa z_2}$,
$z_1\,\frac{\pa}{\pa z_2}+z_2\,\frac{\pa}{\pa z_3}$ and $z_1\,\frac{\pa}{\pa z_3}+z_2\,\frac{\pa}{\pa z_4}$.

\vskip.1in

{\it Case 2.2.1}: $X=z_1\,\frac{\pa}{\pa z_1}-z_2\,\frac{\pa}{\pa z_2}$, or $d\eta=d(z_1\,z_2)\wedge dz_3\wedge dz_4$.
From $i_X\eta=0$ we get $z_1.\,i_{\frac{\pa}{\pa z_1}}\,\eta=z_2.\,i_{\frac{\pa}{\pa z_2}}\,\eta$ $\implies$
there exists a 1-form $\a$, homogeneous of degree one, such that $i_{\frac{\pa}{\pa z_1}}\,\eta=z_2.\,\a$ and $i_{\frac{\pa}{\pa z_2}}\,\eta=z_1.\,\a$.
Note that $i_{\frac{\pa}{\pa z_1}}\,\a=i_{\frac{\pa}{\pa z_2}}\,\a=0$ $\implies$ $\a=A\,dz_3+B\,dz_4$, where $A$ and $B$ are linear.
If we set $\eta=\sum_{1\le i<j\le4}P_{ij}\,dz_i\wedge dz_j$ then
\[
\left\{
\begin{matrix}
z_2.\,(A\,dz_3+B\,dz_4)=i_{\frac{\pa}{\pa z_1}}\,\eta=\sum_{j>1}P_{1j}\,dz_j\\
z_1.\,(A\,dz_3+B\,dz_4)=i_{\frac{\pa}{\pa z_2}}\,\eta=\sum_{j\ne2}P_{2j}\,dz_j\\
\end{matrix}
\right.\,\implies
\]
\[
P_{12}=0\,,\,P_{13}=z_2.\,A\,,\,P_{14}=z_2.\,B\,,\,P_{23}=z_1.\,A\,\text{and}\,P_{24}=z_1.\,B
\]
It follows that
\[
\eta=d(z_1z_2)\wedge(A\,dz_3+B\,dz_4)+C\,dz_3\wedge dz_4\,\,,
\]
where $C=P_{34}$ is quadratic and $A$, $B$ are linear. Now, recall that $L_X\eta=0$ and so
\[
0=L_X\left[d(z_1z_2)\wedge(A\,dz_3+B\,dz_4)+C\,dz_3\wedge dz_4\right]=
\]
\[
=d(z_1z_2)\wedge(X(A)\,dz_3+X(B)\,dz_4)+X(C)\,dz_3\wedge dz_4\,\,\implies
\]
\[
X(A)=X(B)=X(C)=0\,\,\implies
\]
$A=A(z_3,z_4)$, $B=B(z_3,z_4)$ and $C=\wt{C}(z_1.z_2,z_3,z_4)=a.\,z_1.z_2+q(z_3,z_4)$, $q$ quadratic.
Therefore,
\[
\eta=d(z_1\,z_2)\wedge(A(z_3,z_4)\,dz_3+B(z_3,z_4)\,dz_4)+(a.z_1.z_2+q(z_3,z_4))\,dz_3\wedge dz_4\,\,,
\]
In particular, $\eta$ is a pull-back of a 2-form on $\C^3$, $\eta=\Phi^*(\wt\eta)$, where
$\Phi(z_1,z_2,z_3,z_4)=(z_1.z_2,z_3,z_4)=(u,z_3,z_4)$ and
\[
\wt\eta=du\wedge(A(z_3,z_4)\,dz_3+B(z_3,z_4)\,dz_4)+(a\,u+q(z_3,z_4))\,dz_3\wedge dz_4\,\,.
\]
Note also that $\fa_\eta$ and $\fa_{\wt\eta}$ are contained in $\fa_\a$, $\a=A\,dz_3+B\,dz_4$, because $\eta\wedge\a=\wt\eta\wedge\a=0$.

\vskip,1in

{\it Case 2.2.2}: $X=z_1\,\frac{\pa}{\pa z_2}+z_2\,\frac{\pa}{\pa z_3}$, or $d\eta=-z_1\,dz_1\wedge dz_3\wedge dz_4+z_2\,dz_1\wedge dz_2\wedge dz_4$. Here, the integrability relation $i_X\eta=0$ implies that $z_1.\,i_{\frac{\pa}{\pa z_2}}\,\eta=-z_2.\,i_{\frac{\pa}{\pa z_3}}\,\eta$ $\implies$ $i_{\frac{\pa}{\pa z_2}}\,\eta=z_2.\a$ and $i_{\frac{\pa}{\pa z_3}}\,\eta=-z_1.\a$, where $\a=A\,dz_1+B\,dz_4$, $A$ and $B$ linear. With an argument similar to the preceding case, we get
\[
\eta=(z_2\,dz_2-z_1\,dz_3)\wedge(A\,dz_1+B\,dz_4)+C\,dz_1\wedge dz_4\,\,,
\]
where $a\in\C$ and $C$ ishomogeneous of degree two. The reader can check that the condition $L_X\eta=0$ is equivalent to $X(A)=X(B)=0$ and $X(C)+z_2\,B=0$.
Since the first integrals of $X$ are generated by $z_1$, $z_4$ and $z_2^2-2z_1\,z_3$, we get
\[
A=A(z_1,z_4)\,,\,B=B(z_1,z_4)\,\,\text{and}\,\,C=-z_3\,B(z_1,z_4)+a(y^2-2x\,z)+q(z_1,z_4)\,\,,
\]
where $q$ is homogeneous of degree two. As the reader can check, this implies that $\eta=\Phi^*(\wt\eta)$, where
\[
\Phi(z_1,z_2,z_3,z_4)=(z_1,z_4,z^2_2-2z_1\,z_3)=(z_1,z_4,u)\,\,.
\]
and 
\[
\wt\eta=\frac{1}{2}\,du\wedge\left[A\,dz_1+B\,dz_4\right]+\left[a.\,u+q\right]\,dz_1\wedge dz_4
\]

\vskip.1in

{\it Case 2.2.3}: $X=z_1\,\frac{\pa}{\pa z_3}+z_2\,\frac{\pa}{\pa z_4}$. With an argument similar to the preceding cases we get
\[
\eta=(z_2\,dz_3-z_1\,dz_4)\wedge(A\,dz_1+B\,dz_2)+C\,dz_1\wedge dz_2\,\,,
\]
where $A$, $B$ are linear and $C$ homogeneous of degree two. From the condition $L_X\eta=0$ we get
\[
X(A)=X(B)=0\,\,\text{and}\,\,X(C)+z_1\,A+z_2\,B=0\,\,.
\]
Since the first integrals of $X$ are generated by $z_1$, $z_2$ and $z_2\,z_3-z_1\,z_4$ we get $A=A(z_1,z_2)$ $B=B(z_1,z_2)$ and

$C=-z_3\,A(z_1,z_2)-z_4\,B(z_1,z_2)+a(z_2\,z_3-z_1\,z_4)+q(z_1,z_2)$, where $a\in\C$ and $q$ is homogeneous of degree two. Here we obtain
$\eta=\Phi^*(\wt\eta)$ where $\Phi(z_1,z_2,z_3,z_4)=(z_1,z_2,z_2\,z_3-z_1\,z_4)=(z_1,z_2,u)$ and 
\[
\wt\eta=du\wedge\left[A(z_1,z_2)\,dz_1+B(z_1,z_2)\,dz_2\right]+\left[a\,u+q(z_1,z_2)\right]\,dz_1\wedge dz_2\,\,.
\]

\begin{rem}
{\rm We have seen that in the last three cases, 2.2.1, 2.2.2 and 2.2.3, we can choose the coordinates $(x,y,u)\in\C^3$ in such a away that $\eta=\Phi^*(\wt\eta)$, where $\wt\eta=i_Zdx\wedge dy\wedge du$ and $Z$ can be written as
\[
Z(x,y,u)=L_1(x,y)\,\frac{\pa}{\pa x}+L_2(x,y)\,\frac{\pa}{\pa y}+(a.u+q(x,y))\frac{\pa}{\pa u}\,\,,
\]
where $L_1$, $L_2$ are linear and $q$ quadratic. In the generic case, in which the linear part is semi-simple and generic, we can assume that
\[
Z=\la_1\,x\frac{\pa}{\pa x}+\la_2\,y\frac{\pa}{\pa y}+(a\,u+q(x,y))\frac{\pa}{\pa u}\,\,.
\]
If $2\la_1-a\ne0$, $\la_1+\la_2-a\ne0$ and $2\la_2-a\ne0$ then $Z$ is globally linearizable.

Indeed, with these non-resonant conditions there exists a quadratic polynomial $h(x,y)$ such that $Z(h)=a\,h-q$. This implies that
\[
Z\left(u+h(x,y)\right)=a\,u+q(x,y)+Z(h(x,y))=a\,\left(u+h(x,y)\right)\,\,.
\]
In particular, the change of variables $\Phi(x,y,u)=(x,y,u+h(x,y))=(x,y,v)$ linearizes $Z$,
\[
\Phi_*(Z)=\la_1\,x\frac{\pa}{\pa x}+\la_2\,y\frac{\pa}{\pa y}+a\,v\frac{\pa}{\pa v}\,\,.
\]
If we denote the above vector field by $S$ then $\Phi_*\wt\eta=i_S\,dx\wedge dy\wedge dv$. Hence, $\wt{f}=x\,y\,v=x\,y\,(u+h)$ is an integrating factor of $\wt\eta$. It can be chacked that $\eta$ has also an integrating factor.}
\end{rem}

\begin{rem}
{\bf Final remark:} {\rm We would like to observe that the above classification also gives the decomposition into irreducible components of the space of homogeneous foliations of codimension two and degrees zero, one and two.}
\end{rem}



\vskip.2in

\centerline{\Large\bf References}


\vskip.2in

\begin{itemize}

\item[[Ar]] M. Artin: "On the solution of analytic equations"; Inv. Math. 5 (1968), pp. 277-291.

\item[[Br]] M. Brunella: "Sur les feuilletages de l'espace projectif ayant une composante de Kupka"; Enseign. Math. (2) 55 (2009), no. 3-4, 227–234.

\item[[CA-1]]  Calvo Andrade:   "Irreducible components of the space of
holomorphic foliations";  Math. Annalen, no. 299, pp.751-767 (1994).

\item[[CA-2]] O. Calvo-Andrade: "Foliations with a Kupka component on algebraic manifolds"; Bol. Soc. Brasil. Mat. (N.S.) 30 (1999), no. 2, 183–197.

\item[[C-M]] Cano, F.; Mattei, J.-F.: "Hypersurfaces intégrales des feuilletages holomorphes"; Ann. Inst. Fourier (Grenoble) 42 (1992), no. 1-2, 49–72.

\item[[Ce-Cn]] D. Cerveau, F. Cano : "Desingularization of holomorphic foliations and existence of separatrices";  Acta Mathematica, Vol. 169, 1992, pg. 1-103.

\item[[C-LN-S]] C. Camacho, A. Lins Neto e P. Sad: "Foliations with algebraic limit sets";
Ann. of Math. 136 (1992), pg. 429-446.

\item[[Ce-LN]] D. Cerveau, A. Lins Neto: "Irreducible
components of the space of holomorphic foliations of degree two in $\mathbb C P(n)$, $n \ge 3$";  Ann. of Math. (1996) 577-612.

\item[[Ce-LN 1]] D. Cerveau, A. Lins Neto: "Codimension-one foliations in $\mathbb C P(n)$ $n
\ge 3$ with Kupka components"; Ast\'erisque 222 (1994) pg. 93-132.

\item[[Ce-LN-2]] D. Cerveau \& A. Lins Neto: "A structural theorem for codimension-one foliations on $\p^n$, $n\ge3$, with an application to degree three foliations"; Ann. Sc. Norm. Super. Pisa Cl. Sci. (5) 12 (2013), no. 1, 1–41.

\item[[Ce-LN-CA-G]] D. Cerveau, A. Lins Neto, J. O. Calvo Andrade, L. Giraldo : "Irreducible components of the space of foliations associated to the affine Lie algebra"; To appear in Ergodic Theory and Dynamical Systems.

\item[[Ce-LN-Lo-Pe-Tou]] D. Cerveau, A. Lins Neto, F. Loray, J. V. Pereira, F. Touzet: "Complex codimension one foliations and Godbillon-Vey sequences"; Mosc. Math. J. 7 (2007), no. 1, 21–54, 166. 

\item[[Ce-Ma]] D. Cerveau, J.-F. Mattei: "Formes int\'egrables holomorphes
Singuli\`eres"; Ast\'erisque, vol.97 (1982).

\item[[DR]] G. De Rham : "Sur la division des formes et des courants par une forme lin\'eaire";
Comm. Math. Helvetici, 28 (1954), pp. 346-352.

\item[[J]]  J.P. Jouanolou : "\'Equations de
Pfaff alg\`ebriques"; Lecture Notes in
Math. 708, Springer-Verlag, Berlin, 1979.

\item[[K]] I. Kupka: "The singularities of integrable structurally stable
Pfaffian forms"; Proc. Nat. Acad. Sci. U.S.A., 52 (1964), pg. 1431-1432.

\item[[K-N]] Dan, Krishanu \& Nagaraj, D. S.: "Null correlation bundle on projective three space"; J. Ramanujan Math. Soc. 28A (2013), 75-80.

\item[[LN]] A. Lins Neto: "Germs of complex two dimensional foliations", Accepted for publication in the "Bulletin of the Brazilian Mathematical Society".

\item[[LN 1]] A. Lins Neto: "Algebraic solutions of polynomial differential equations and foliations in dimension two."; Holomorphic dynamics (Mexico, 1986), 192--232, Lecture Notes in Math., 1345, Springer, Berlin, 1988.

\item[[LN-So]] A. Lins Neto, M.G. Soares: "Algebraic solutions of
onde-dimensional foliations"; J. Diff. Geometry 43 (1996) pg. 652-673.

\item[[L-P-T]] F. Loray, J. V. Pereira, F. Touzet: "Foliations with trivial canonical bundle on Fano 3-folds"; Math. Nachr. 286 (2013), no. 8-9, 921–940. 

\item[[M-1]] B. Malgrange : "Frobenius avec singularit\'es I. Codimension un."
Publ. Math. IHES, 46 (1976), pp. 163-173.

\item[[M-2]] B. Malgrange: "Frobenius avec singularités. II. Le cas général"; Invent. Math. 39 (1977), no. 1, 67–89.

\item[[M-M-1]] Marín, D.; Mattei, Jean-François: "Topology of singular holomorphic foliations along a compact divisor"; J. Singul. 9 (2014), 122–150.

\item[[M-M-2]] Marín, D.; Mattei, Jean-François: "Monodromy and topological classification of germs of holomorphic foliations"; Ann. Sci. Éc. Norm. Supér. (4) 45 (2012), no. 3, 405–445.

\item[[Ma]] Mattei, J.-F.: "Modules de feuilletages holomorphes singuliers. I. Équisingularité"'. Invent. Math. 103 (1991), no. 2, 297–325.

\item[[Ma-Mo]] J.F. Mattei \& R. Moussu: "Holonomie et int\'egrales premi\`eres";
Ann. Ec. Norm. Sup. 13 (1980), pg. 469-523. 

\item[[Me]] Medeiros, Airton S.: "Singular foliations and differential p-forms"; Ann. Fac. Sci. Toulouse Math. (6) 9 (2000), no. 3, 451–466.

\item[[Se]] A. Seidenberg: "Reduction of singularities of the differential
equation $Ady=Bdx$"; Amer. J. de Math. 90 (1968), pg. 248-269.

\item[[Sa]] K. Saito: "On a generalization of de Rham lemma"; Annales de l'institut Fourier, 26 no. 2 (1976), p. 165-170

\end{itemize}


\vskip.3in

{\sc D. Cerveau}

{\sl Inst. Math\'ematique de Rennes}

{\sl Campus de Beaulieu}

{\sl 35042 RENNES Cedex}

{\sl Rennes, France}

{\tt E-Mail: dominique.cerveau@univ-rennes1.fr}

\vskip.3in

{\sc A. Lins Neto}

{\sl Instituto de Matem\'atica Pura e Aplicada}

{\sl Estrada Dona Castorina, 110}

{\sl Horto, Rio de Janeiro, Brasil}

{\tt E-Mail: alcides@impa.br}

\end{document}